\documentclass[a4paper,10pt]{amsart}
\usepackage{amssymb}
\usepackage{psfig}

\input xypic
\xyoption{all}
\xyoption{arc}

\def\sursection#1{\refstepcounter{section}
\vskip0.6cm\begin{centerline}{\bf \large \arabic{section}. #1}
\end{centerline}\vskip0.6cm} 

\newcounter{soussection}[section]

\def\soussection#1{\refstepcounter{soussection}
\noindent{\bf \arabic{section}.\Alph{soussection}. #1.}} 
\newcounter{soussoussection}[soussection]

\def\soussoussection#1{\refstepcounter{soussoussection}
\noindent{\it \arabic{section}.\Alph{soussection}.\arabic{soussoussection}. #1.}} 

\newtheorem{theo}{Theorem}[section]

\newtheorem{prop}[theo]{Proposition}

\newtheorem{lem}[theo]{Lemma}

\newtheorem{coro}[theo]{Corollary}

\def\remark#1{{\refstepcounter{theo}\label{#1}\noindent\sc Remark  
\arabic{section}.\arabic{theo} - }}
\def\example#1{{\refstepcounter{theo}\label{#1}\noindent\sc Example 
\arabic{section}.\arabic{theo} - }}

\def\equat{\refstepcounter{theo}$$~}
\def\endequat{\leqno{\boldsymbol{(\arabic{section}.\arabic{theo})}}~$$}

\newcounter{numero}[section]

\def\FM{{\mathbb{F}}}

\def\QM{{\mathbb{Q}}}

\def\ZM{{\mathbb{Z}}}

\def\SG{{\mathfrak S}}

\def\a{\alpha}
\def\b{\beta}

\def\g{\gamma}

\def\d{\delta}
\def\D{\Delta}
\def\e{\varepsilon}

\def\ch{\chi}
\def\l{\lambda}

\def\m{\mu}

\def\o{\omega}
\def\O{\Omega}

\def\s{\sigma}

\def\t{\tau}

\def\z{\zeta}

\def\AC{{\mathcal{A}}}

\def\CC{{\mathcal{C}}}

\def\QC{{\mathcal{Q}}}

\def\Ab{{\mathbf A}}
\def\Bb{{\mathbf B}}

\def\Db{{\mathbf D}}

\def\Gb{{\mathbf G}}

\def\Lb{{\mathbf L}}

\def\Ob{{\mathbf O}}
\def\Pb{{\mathbf P}}

\def\Sb{{\mathbf S}}
\def\Tb{{\mathbf T}}

\def\Xb{{\mathbf X}}

\def\Zb{{\mathbf Z}}

\def\pb{{\mathbf p}}

\def\gti{{\tilde{g}}}

\def\sti{{\tilde{s}}}
\def\tti{{\tilde{t}}}

\def\xti{{\tilde{x}}}
\def\yti{{\tilde{y}}}

\def\sba{{\bar{s}}}
\def\tba{{\bar{t}}}

\def\Bbt{{\tilde{\Bb}}}

\def\Gbt{{\tilde{\Gb}}}

\def\Tbt{{\tilde{\Tb}}}

\def\alpt{{\tilde{\alpha}}}

\def\Delt{{\tilde{\Delta}}}

\def\Phit{{\tilde{\Phi}}}

\def\omet{{\tilde{\omega}}}

\def\Aut{\mathop{\mathrm{Aut}}\nolimits}

\def\diag{\mathop{\mathrm{diag}}\nolimits}

\def\Id{\mathop{\mathrm{Id}}\nolimits}

\def\Im{\mathop{\mathrm{Im}}\nolimits}

\def\Ker{\mathop{\mathrm{Ker}}\nolimits}

\def\mod{\mathop{\mathrm{mod}}\nolimits}

\def\Oplus{\mathop{\oplus}}

\def\tete#1{\par\leavevmode\makebox[0.7cm]{$(\mathrm{#1})$}}

\def\imp{\Rightarrow}

\def\to{\rightarrow}
\def\longto{\longrightarrow}
\def\injto{\hookrightarrow}

\def\longmapright#1{\hspace{0.3em}\smash{
     \mathop{\longrightarrow}\limits^{#1}}\hspace{0.3em}}

\def\longmapleftright#1{\hspace{0.3em}\smash{
     \mathop{\longleftrightarrow}\limits^{#1}}\hspace{0.3em}}

\def\fonction#1#2#3#4#5{\begin{array}{rccc}
{#1} : & {#2} & \longto & {#3} \\
& {#4} & \longmapsto & {#5} 
\end{array}}

\def\incl{\hspace{0.05cm}{\subset}\hspace{0.05cm}}

\def\vide{\varnothing}

\def\DS{\displaystyle}
\def\SS{\scriptstyle}

\def\fin{~$\SS \blacksquare$}
\def\finl{~$\SS \square$}

\def\lexp#1#2{\kern\scriptspace\vphantom{#2}^{#1}\kern-\scriptspace#2}
\def\le{\hspace{0.1em}\mathop{\leqslant}\nolimits\hspace{0.1em}}
\def\ge{\hspace{0.1em}\mathop{\geqslant}\nolimits\hspace{0.1em}}

\mathchardef\lllllll="3278
\def\SEC{$\lllllll$}
\mathchardef\inferieur="321E
\mathchardef\superieur="321F

\def\med{\medskip}

\def\eqna{\begin{eqnarray*}}
\def\endeqna{\end{eqnarray*}}

\def\proof{\noindent{\sc{Proof}~-} }

\def\torus{maximal torus }

\def\borel{Borel subgroup }

\def\para{parabolic subgroup }

\def\levi{Levi subgroup }

\def\resp{respectively }
\def\ssi{if and only if }

\def\mini{{\mathrm{min}}}
\def\itemth#1{\item[${\mathrm{(#1)}}$]}

\def\Div{{\mathrm{Div}}}
\def\itemth#1{\item[${\mathrm{(#1)}}$]}

\def\simp{{\mathrm{sc}}}
\def\adj{{\mathrm{ad}}}
\def\ZMP{{\ZM_{(p)}}}

\def\vertical{\vphantom{\frac{\DS{A^A_A}}{\DS{A_A^A}}}}
\def\Com{\mathop{\mathrm{Com}}\nolimits}

\begin{document}

\title{Quasi-isolated elements in reductive groups}

\author{C\'edric Bonnaf\'e}

\address{C\'edric Bonnaf\'e, 
CNRS - UMR 6623, Universit\'e de Franche-Comt\'e, 16 Route de Gray, 25030 
BESAN\c{C}ON Cedex, FRANCE}

\email{bonnafe@math.univ-fcomte.fr}

\date{\today}

\subjclass{20}

\markboth{C. Bonnaf\'e}{Quasi-isolated elements}

\bigskip

\bigskip

\begin{abstract}
A semisimple element $s$ of a connected reductive group $\Gb$ is said 
{\it quasi-isolated} (respectively {\it isolated}) if $C_\Gb(s)$ 
(respectively $C_\Gb^\circ(s)$) is not contained in a \levi of a proper 
\para of $\Gb$. We study properties of quasi-isolated semisimple elements and 
give a classification in terms of the affine Dynkin diagram of $\Gb$. 
Tables are provided for adjoint simple groups.
\end{abstract}

\maketitle

\bigskip

\bigskip

\sursection{Preliminaries and notation}

\bigskip

\soussection{Notation} 
Let $\FM$ be an algebraically closed field. Let $p$ denote its characteristic. 
By a variety (respectively an algebraic group), we mean an algebraic variety 
(respectively an algebraic group) over $\FM$. We denote by $\ZMP$ the localization 
of $\ZM$ at the prime ideal $p\ZM$ (in particular, if $p=0$, then $\ZMP=\QM$). 

We fix a connected reductive group $\Gb$. We denote by $\Zb(\Gb)$ its center and 
$\Db(\Gb)$ its derived subgroup. If $g \in \Gb$, we denote by $g_s$ 
(respectively $g_u$) its semisimple (respectively unipotent) part, 
$C_\Gb(g)$ its centralizer and $C_\Gb^\circ(g)$ the neutral component 
of $C_\Gb(g)$. We denote by $o(g) \in \{1,2,3,\dots\} \cup \{\infty\}$ 
the order of $g$.

\bigskip

\soussection{Isolated and quasi-isolated elements} 
The element $g \in \Gb$ is said {\it quasi-isolated} 
(respectively {\it isolated}) if $C_\Gb(g_s)$ (respectively 
$C_\Gb^\circ(g_s)$) is not contained in a \levi of a proper 
\para of $\Gb$. If there is some ambiguity, we will speak about 
$\Gb$-isolated or $\Gb$-quasi-isolated elements to refer to the ambient group. 
Of course, an isolated element is quasi-isolated.

The isolated elements are present in many different papers while 
the quasi-isolated ones are not often mentionned (see \cite[\SEC 4.5]{bonnafe mackey}). 
One reason might be the following~: if the derived group of $\Gb$ is 
simply connected, then centralizers of semisimple elements are connected 
(by a theorem of Steinberg \cite[Theorem 8.1]{steinberg}, see 
also \cite[Chapter VI, \SEC 2, Exercise 1]{bourbaki}) 
so the notions of isolated and quasi-isolated elements coincide. 
Another possible reason is that the notion of isolated element depends 
only on the Dynkin diagram of $\Gb$, by opposition to the notion of 
quasi-isolated element (see Proposition \ref{iso isole} and Example 
\ref{typique 2}).

Whenever the derived group of $\Gb$ is not simply connected, a quasi-isolated 
element might not be isolated. The following extreme case can even happen~: 
there exist quasi-isolated semisimple elements $s$ which are {\it regular} 
(that is, such that $C_\Gb^\circ(s)$ is a maximal torus), as it is shown 
by the following example. 

\bigskip

\example{typique} 
Let $n \ge 2$ be a natural number invertible in $\FM$. 
Let us assume in this example that $\Gb=\Pb\Gb\Lb_n(\FM)$. Let $\z$ 
be a primitive $n$-th root of unity in $\FM$ and let $s$ be the image 
of $\diag(1,\z,\z^2,\dots,\z^{n-1}) \in\Gb\Lb_n(\FM)$ in $\Gb$. 
Then $C_\Gb^\circ(s)$ is the maximal torus consisting of the image 
of diagonal matrices in $\Gb$ (in particular, $s$ is regular, so it is 
not isolated) but $C_\Gb(s)/C_\Gb^\circ(s)$ is cyclic of order $n$~: 
it is generated by a Coxeter element of the Weyl group of $\Gb$ relatively to 
$C_\Gb^\circ(s)$. Therefore, $s$ is quasi-isolated.\finl

\bigskip

\soussection{Root system} 
The notions of isolated and quasi-isolated elements involve only the 
semisimple part, so we will focus on semisimple elements. For this reason, 
we fix once and for all a maximal torus of $\Gb$~: determining if an element 
of this torus is quasi-isolated or not can be done thanks to the root 
system or the Weyl group relatively to this torus.

\medskip

Let $\Bb$ be a \borel of $\Gb$ and let $\Tb$ be a \torus of $\Bb$. 
Let $W$ be the Weyl group and let $\Phi$ be the root system of $\Gb$ 
relatively to $\Tb$. Let $\Phi^+$ (respectively $\D$) denote the positive 
root system (respectively the basis) of $\Phi$ associated to $\Bb$. 

We fix once and for all an element $s \in \Tb$. We denote by $\Phi(s)$ and by 
$W^\circ(s)$ respectively the root system and the Weyl group 
of $C_\Gb^\circ(s)$ relatively to $\Tb$. We set~:
$$W(s)=\{w \in W~|~\lexp{w}{s}=s\}.$$
Let $\Bb(s)$ be a \borel of $C_\Gb^\circ(s)$ containing $\Tb$ 
and let $\Phi^+(s)$ (respectively $\D(s)$) denote the positive 
root system (respectively the basis) of $\Phi(s)$ 
associated to $\Bb(s)$. We set~:
$$A(s)=\{w \in W(s)~|~w(\Phi^+(s))=\Phi^+(s)\}.$$

\bigskip

\example{borel cgos} $C_\Bb^\circ(s)$ is a \borel of $C_\Gb^\circ(s)$ 
containing $\Tb$. If $\Bb(s)=C_\Bb^\circ(s)$, then 
$\Phi^+(s)=\Phi^+ \cap \Phi(s)$.\finl

\bigskip

We gather some elementary facts~:

\bigskip

\begin{prop}\label{collection}
Let $s \in \Tb$. Then~:
\begin{itemize}
\itemth{a} $\Phi(s)=\{\a \in \Phi~|~\a(s)=1\}$.

\itemth{b} $W(s)$ is the Weyl group of $C_\Gb(s)$ relatively to $\Tb$.

\itemth{c} $W(s) = A(s) \ltimes W^\circ(s)$.

\itemth{d} $A(s) \simeq C_\Gb(s)/C_\Gb^\circ(s)$.
\end{itemize}
\end{prop}

\bigskip

\begin{coro}\label{isolitude}
Let $s \in \Tb$. Then~:
\begin{itemize}
\itemth{a} $s$ is isolated (respectively quasi-isolated) 
\ssi $W^\circ(s)$ (respectively $W(s)$) is not contained in 
a proper \para of $W$.

\itemth{b} The following are equivalent~:
\begin{itemize}
\itemth{1} $s$ is isolated~;

\itemth{2} $\Phi(s)$ is not contained in a proper parabolic subsystem of $\Phi$~;

\itemth{c} $|\D(s)|=|\D|$. 
\end{itemize}
\end{itemize}
\end{coro}

\bigskip

\begin{prop}\label{ordre fini pareil}
Let $s \in \Tb$. Then there exists an element $s' \in \Tb$, 
of finite order, such that $C_\Gb(s)=C_\Gb(s')$.
\end{prop}

\bigskip

\proof Let $\Sb$ denote the Zarisky closure of the group generated 
by $s$. Then $\Sb/\Sb^\circ$ is generated by the image of $s$, so it 
is cyclic. Moreover, $C_\Gb(s)=C_\Gb(\Sb)$. Therefore, Proposition 
\ref{ordre fini pareil} follows immediately from the following easy lemma~:

\medskip

\begin{lem}\label{T V}
Let $\Db$ be a diagonalizable group acting on an affine variety $\Xb$. 
We assume that $\Db/\Db^\circ$ is cyclic. Then there exists an element 
$t \in \Db$ of finite order such that $\Xb^\Db=\Xb^t$.
\end{lem}

\bigskip

\noindent{\sc Proof of Lemma \ref{T V} - } We first prove the following statement~: 

\bigskip

\begin{quotation}
{\small \noindent{$\boldsymbol{(*)}$}\quad{\it For 
every prime number $\ell\not= p$, there exists an element 
$t \in \Db^\circ$ of $\ell$-power order such that $\Xb^{\Db^\circ}=\Xb^t$.}

\medskip

\noindent{\sc Proof of $(*)$ - } By \cite[Proposition 1.12]{borel}, there 
exists a finite dimensional rational representation $V$ of $\Db^\circ$ and 
a $\Db^\circ$-equivariant closed embedding $\Xb \injto V$. So 
$\Xb^{\Db^\circ} = V^{\Db^\circ} \cap \Xb$. If $\ch \in X(\Db^\circ)$, we set
$$V_\ch=\{v \in V~|~\forall t \in \Db,~t.v=\ch(t)v\}.$$
Let $\ch_1$,\dots, $\ch_k \in X(\Db^\circ)$ be the distinct 
non-zero weights of $\Db^\circ$ in its action on $V$. Then~:
$$V=V^{\Db^\circ} \oplus \Bigl(\Oplus_{i=1}^k V_{\ch_i}\Bigr).$$
Now, the subgroup $L$ of $\Db^\circ$ consisting of elements 
of $\ell$-power order is Zarisky dense in $\Db^\circ$. Therefore, there 
exists $t \in L$ such that $t \not\in \Ker \ch_1 \cup \dots \cup \Ker \ch_k$. 
This implies that $V^{\Db^\circ}=V^t$. This shows $(*)$.\finl}
\end{quotation}

\bigskip

\noindent Now let $D$ be a finite cyclic 
subgroup of $\Db$ such that $\Db = D \times \Db^\circ$. Let $d \in D$ be such 
that $D=<d>$. Then $\Xb^\Db=(\Xb^{\Db^\circ})^d$. Let $\ell$ be a prime 
number greater than $|D|$. Then, by $(*)$, there exists $t \in \Db^\circ$ 
of $\ell$-power order such that $\Xb^{\Db^\circ}=\Xb^t$. Then, since 
$t$ and $d$ have coprime order, we have 
$\Xb^{td}=(\Xb^t)^d=(\Xb^{\Db^\circ})^d=\Xb^\Db$.\fin

\bigskip

\noindent{\sc Remark - } 
Lemma \ref{T V} and statement $(*)$ are slight refinements 
of a well-known lemma on the action 
of tori on affine varieties (see for instance \cite[Proposition 0.7]{dmbook}). 
Moreover, Lemma \ref{T V} does not remain valid if $\Db/\Db^\circ$ is not 
cyclic. For instance, assume that $p \not=2$ and 
consider the action of $\Db=\{1,-1\} \times \{1,-1\}$ 
on $\Ab^3(\FM)$ by $(\e,\e').(x,y,z)=(\e x, \e' y, \e\e' z)$.\finl

\bigskip

\sursection{Isotypic morphisms}

\bigskip

\soussection{Definition\label{def}}
A morphism $\pi : \Gbt \to \Gb$ is said {\it isotypic} if $\Gbt$ is a connected 
reductive group, if $\Ker \pi$ is central in $\Gbt$ and if $\Im \pi$ 
contains the derived group of $\Gb$. 

\bigskip

\noindent{\sc Example and notation - } Let $\pi_\simp : \Gb_\simp \to \Gb$ 
be a simply connected covering of the derived group of $\Gb$. Let $\Gb_\adj$ 
denote the adjoint group of $\Gb$ and let $\pi_\adj : \Gb \to \Gb_\adj$ 
be the canonical 
surjective morphism. Then $\pi_\simp$ and $\pi_\adj$ are 
isotypic morphisms. We set 
$\Bb_\adj=\pi_\adj(\Bb)$ and $\Tb_\adj=\pi(\Tb)$. Then $\Bb_\adj$ is a 
\borel of $\Gb_\adj$ and $\Tb_\adj$ is a \torus of $\Bb_\adj$. 
Moreover, if $t \in \Tb$, we set $\tba=\pi_\adj(t) \in \Tb_\adj$.\finl 

\medskip

We fix in this section an isotypic morphism $\pi : \Gbt \to \Gb$.  
Let $\Ker' \pi = \Db(\Gbt) \cap \Ker \pi$. It must be noticed that 
$\Ker' \pi$ is a finite abelian group of order prime to $p$. 
Let $\Bbt=\pi^{-1}(\Bb)$ 
and $\Tbt=\pi^{-1}(\Tb)$. Then $\Bbt$ is a \borel of $\Gbt$ and $\Tbt$ is a \torus 
of $\Bbt$. We will identify the Weyl group of $\Gbt$ relatively to $\Tbt$ with $W$ 
through the morphism $\pi$. Let $\Phit$ denote the root system of $\Gbt$ relatively to 
$\Tbt$. Then the morphism $\pi^* : X(\Tb) \to X(\Tbt)$ induced by $\pi$ 
provides a bijection $\Phi \longmapleftright{\sim} \Phit$. 

\medskip

Since $\Im \pi$ contains $\Db(\Gb)$, we have $\pi(\Gbt).\Zb(\Gb)^\circ=\Gb$. 
We fix once and for all in this section an element 
$\sti \in \Tbt$ such that $\pi(\sti) \in s\Zb(\Gb)$. Then 
\equat\label{pi cg}
\pi(C_\Gbt(\sti)).\Zb(\Gb)^\circ \incl C_\Gb(s).
\endequat
Moreover, by Proposition \ref{collection} (a), we have 
\equat\label{pi cgo}
\pi(C_\Gbt^\circ(\sti)).\Zb(\Gb)^\circ= C_\Gb^\circ(s).
\endequat
Therefore, $W(\sti) \incl W(s)$ and $W^\circ(\sti)=W^\circ(s)$. 
Moreover, $A(\sti) \incl A(s)$ (if we choose $\Bbt(\sti)=\pi^{-1}(\Bb(s))$).
These remarks have the following consequences~:

\bigskip

\begin{prop}\label{iso isole}
With the above notation, we have~:
\begin{itemize}
\itemth{a} If $\sti$ is quasi-isolated in $\Gbt$, then $s$ is quasi-isolated in $\Gb$.

\itemth{b} $\sti$ is isolated in $\Gbt$ \ssi $s$ is isolated in $\Gb$.
\end{itemize}
\end{prop}

\bigskip

The following example shows that the converse to statement (a) of Proposition 
\ref{iso isole} is not true in general.

\bigskip

\example{typique 2} Keep here the hypothesis and notation of Example \ref{typique}. 
Assume that $\Gbt=\Gb\Lb_n(\FM)$ and that $\pi : \Gbt \to \Gb$ is the 
canonical morphism. Let $\sti=\diag(1,\z,\z^2,\dots,\z^{n-1})$. 
Then $\sti$ is not quasi-isolated in $\Gbt$ since $C_\Gbt(\sti)$ is 
a maximal torus. But $s=\pi(\sti)$ is quasi-isolated in $\Gb$ 
as it is shown in Example \ref{typique}.\finl

\bigskip

\remark{injectif pi} 
If $\pi$ is injective, then the inclusion \ref{pi cg} is an equality.  
So $s$ is quasi-isolated in $\Gb$ \ssi $\sti$ is quasi-isolated in $\Gbt$.\finl

\bigskip

\soussection{The groups ${\boldsymbol{A(s)}}$ et ${\boldsymbol{A(\sti)}}$} 
We will compare here the groups $A(s)$ and $A(\sti)$ 
in order to obtain general properties of the group $A(s)$. 
Most of the results of this subsection are well-known, particularly 
the Corollary \ref{injection plus}  
(see \cite[lemme 9.2]{steinberg} and \cite[lemme 2.1]{broue michel}) 
but they are rarely stated in the whole generality of this subsection. 

Let $\Com(\Gb)$ denote the set of couples $(x,y) \in \Gb \times \Gb$ such that 
$xy=yx$. This is a closed subvariety of $\Gb \times \Gb$. 
If $(x,y) \in \Com(\Gb)$, we denote by $\o(x,y)$ the element 
$[\xti,\yti]=\xti\yti\xti^{-1}\yti^{-1} \in \Gbt$ where $\xti \in \Gbt$ and 
$\yti$ are two elements of $\Gbt$ such that 
$\pi(\xti) \in x \Zb(\Gb)$ and $\pi(\yti) \in y \Zb(\Gb)$. 
It is easily checked that 
$\o(x,y)$ depends only on $x$ and $y$ and does not depend on the choice of 
$\xti$ and $\yti$. Moreover, $\pi([\xti,\yti])=[x,y]=1$ so 
$\o(x,y) \in \Ker' \pi$.

\bigskip

\begin{lem}\label{bilineaire}
Let $x$, $x'$, $y$ and $y'$ be four elements of $\Gb$ such that $xy=yx$, 
$x'y=yx'$ and $xy'=y'x$. Then~:
$$\o(x',yy')=\o(x,y)\o(x,y'),$$
$$\o(xx',y)=\o(x,y)\o(x',y)$$
$$\o(x,y)=\o(y,x)^{-1}.\leqno{\mathit{and}}$$
\end{lem}

\bigskip

\proof 
Let us show the first equality (the second can be shown similarly and the 
third one is obvious). 
Let $\xti$, $\yti$ and $\yti'$ be three elements of $\Gbt$ such that 
$\pi(\xti) \in x\Zb(\Gb)$, $\pi(\yti) \in y\Zb(\Gb)$ and $\pi(\yti') \in y'\Zb(\Gb)$. 
Then $\pi(\yti\yti') \in yy' \Zb(\Gb)$ and so 
\eqna
\o(x,yy')&=& [\xti,\yti\yti'] \\
                &=& \xti\yti\yti'\xti^{-1}\yti^{\prime -1}\yti^{-1} \\
		&=& \xti\yti\xti^{-1}\xti\yti'\xti^{-1}\yti^{\prime -1}\yti^{-1} \\
		&=& \xti\yti\xti^{-1}\o(x,y')\yti^{-1} \\
		&=& \o(x,y) \o(x,y'),
\endeqna
where the last equality follows from the fact that $\o(x,y')$ is 
central in $\Gbt$.\fin

\bigskip

Let $\o_s : C_\Gb(s) \to \Ker' \pi$, $g \mapsto \o(g,s)$. 
The Lemma \ref{bilineaire} shows that $\o_s$ is a morphism of groups.  

\bigskip

\begin{lem}\label{noyau}
$\Ker \o_s = \pi(C_\Gbt(\sti)).\Zb(\Gb)^\circ$.
\end{lem}

\bigskip

\proof Let $g \in \Ker \o_s$. There exists $\gti \in \Gbt$ such that 
$\pi(\gti) \in g \Zb(\Gb)^\circ$. Since $\o_s(g)=[\gti,\sti]=1$, we have 
$\gti \in C_\Gbt(\sti)$.\fin

\bigskip

\begin{coro}\label{injection}
We have~:
\begin{itemize}
\itemth{a} $\o_s$ induces a morphism $\omet_s : A(s) \longto \Ker' \pi$. 
We have $\Ker \omet_s = A(\sti)$ and 
$\Im \omet_s =\Im \o_s= \{z \in \Ker \pi~|~\sti$ and $\sti z$ 
are conjugated in $\Gbt\}$.

\itemth{b} $|A(s)/A(\sti)|$ is a finite abelian group of order dividing 
$|\Ker' \pi|$ (so prime to $p$).
\end{itemize}
\end{coro} 

\bigskip

\proof (a) By Lemma \ref{noyau} and equality \ref{pi cgo}, 
$C_\Gb^\circ(s)$ is contained in the kernel of $\o_s$. 
Since $A(s) \simeq C_\Gb(s)/C_\Gb^\circ(s)$, we get the first assertion. 
The second follows again from Lemma 
\ref{noyau}. 

Let us show the last one. Let $A=\{z\in \Ker \pi~|~\sti$ and $\sti z$ are conjugated 
in $\Gbt\}$. Let $g \in C_\Gb(s)$. Let $\gti$ be an element of $\Gbt$ such that 
$\pi(\gti) \in g \Zb(\Gb)$. Set 
$z=\o_s(g)$. Then $\sti z = \gti \sti \gti^{-1}$, which shows that 
$z \in A$. So $\Im \o_s \incl A$. 
Conversely, let $z \in A$. Then there exists $\gti \in \Gbt$ such that 
$\sti z = \gti \sti \gti^{-1}$. Set $g=\pi(\gti)$. Then 
$z=[\gti,\sti]$, which shows that $g \in C_\Gb(s)$ and that $z=\o_s(g)$. 
So $A \incl \Im \o_s$.

\medskip

\tete{b} follows immediately from (a).\fin

\bigskip

\begin{coro}\label{injection plus}
The group $A(s)$ is isomorphic to a subgroup of the $p'$-part of the 
fundamental group of $\Db(\Gb)$. The exponent of $A(s)$ divides 
the order of $\sba$ in $\Gb_\adj$, whenever this one is finite.
\end{coro}

\bigskip

\proof This statement does not involve the group $\Gbt$. So we can choose 
for $\pi$ the most convenient morphism for this question. We thus assume 
that $\pi : \Gbt \to \Gb$ is the morphism $\pi_\simp : \Gb_\simp \to \Gb$ 
defined in \SEC\ref{def}. Then $\Ker' \pi=\Ker \pi$ is the $p'$-part 
of the fundamental group of $\Gb$. 

Moreover, Steinberg's Theorem \cite[th\'eor\`eme 8.1]{steinberg} tells 
us that $C_\Gbt(\sti)$ is connected, so $A(\sti)=1$. 
So the first assertion follows immediately from Corollary \ref{injection} (b). 
Let us show now the second assertion. Let $n$ denote the order of $\sba$ in 
$\Gb_\adj$ and let $g \in C_\Gb(s)$. We must show that $g^n \in C_\Gb^\circ(s)$. 
Then, by Lemma \ref{bilineaire}, 
we have $\o_s(g^n)=\o(g,s)^n=\o(g,s^n)$. But, $\sti^n$ is central 
in $\Gbt$. Thus $\o_s(g^n)=1$, which shows that 
$g^n \in \pi(C_\Gbt(\sti))\Zb(\Gb)^\circ$ by Lemma \ref{noyau}. 
But, again by Steinberg's Theorem,
we have $C_\Gbt(\sti)=C_\Gbt^\circ(\sti)$, so $g^n \in C_\Gb^\circ(s)$ 
by \ref{pi cgo}.\fin

\bigskip

\soussection{Isotypic morphisms and quasi-isolated elements} 
The Proposition \ref{iso isole} shows that the notion of isolated element 
depends only on the isogeny class $\Gb$. 
On the other hand, the Example \ref{typique 2} shows that 
the notion of quasi-isolated element does not behave so nicely. 
We will use the morphism $\o_s$ to study a weak converse to the statement 
(a) of Proposition \ref{iso isole}. 
This weak converse will also be used to obtain some classification result 
for quasi-isolated elements. 

Let $e_s^\pi$ denote the exponent of the group $A(s)/A(\sti)$ 
(recall that $e_s^\pi$ divides the exponent of $\Ker' \pi$ and the order 
of $\sba$ in $\Gb_\adj$). A result analogous to the following has been shown in 
\cite[preuve du corollaire 4.5.3]{bonnafe mackey}.

\bigskip

\begin{prop}\label{centralisateur}
The group $C_\Gb(s)$ is contained in $\pi(C_\Gbt(\sti^{e_s^\pi})).\Zb(\Gb)^\circ$.
\end{prop}

\bigskip

\proof Let $g \in C_\Gb(s)$. Then $\o_s(g)^{e_s^\pi}=1$. But, by Lemma 
\ref{bilineaire}, we have $\o_s(g)^{e_s^\pi}=\o_{s^{e_s^\pi}}(g)$. 
This shows that $g \in \Ker \o_{s^{e_s^\pi}} = 
\pi(C_\Gbt(\sti^{e_s^\pi})).\Zb(\Gb)^\circ$ (see Lemma \ref{noyau}).\fin

\bigskip

\begin{coro}\label{quasi isole}
If $s$ is quasi-isolated in $\Gb$, then $\sti^{e_s^\pi}$ is quasi-isolated in $\Gbt$.
\end{coro}

\bigskip

\begin{coro}\label{iso qiso}
Let $e$ be the exponent of $\Ker \pi_\simp$. If $s$ is quasi-isolated in $\Gb$, then 
$s^e$ is isolated in $\Gb$.
\end{coro}

\bigskip

\proof Once again, the group $\Gbt$ is not involved in this statement, 
so we can assume here that $\pi=\pi_\simp$. Then $e_s^\pi$ divides $e$ so, 
by Corollary \ref{quasi isole}, $\sti^e$ is quasi-isolated 
in $\Gbt=\Gb_\simp$. But, since $\Gbt$ is simply connected, $\sti^e$ 
is isolated in $\Gbt$. Therefore, by Proposition \ref{isolitude} (a), 
$s^e$ is isolated in $\Gb$.\fin

\bigskip

\begin{coro}\label{ordre fini}
If $s$ is quasi-isolated, then $\sba$ has finite order.
\end{coro}

\bigskip

\proof By Corollary \ref{iso qiso}, 
we may assume that $s$ is isolated in $\Gb$. Let  
$\m : X(\Tb_\adj) \to \FM^\times$, $\chi \mapsto \chi(\sba)$. 
Then $\m$ is a morphism of groups and $\Phi(s)=\Phi \cap \Ker \m$ 
(here, we identify $X(\Tb_\adj)$ to a subgroup of $X(\Tb)$ via 
the morphism $\pi_\adj$). 
Since $s$ is isolated, it follows from Corollary \ref{isolitude} (b) 
that $<\Phi(s)>$ has finite index in $X(\Tb_\adj)$. 
So $\Ker \m$ has finite index in $X(\Tb_\adj)$. Let $d$ 
denote this index (that is the order of the image of $\m$). 
Then $\sba$ has order $d$ in $\Tb_\adj$.\fin 

\bigskip

\sursection{Semisimple elements of finite order}

\bigskip

We will describe in this subsection the possible structure of 
the centralizer of a semisimple element in $\Gb$. By Proposition 
\ref{ordre fini pareil}, we can focus on semisimple elements 
of finite order. 
For this, we fix an injective morphism $\imath : (\QM/\ZM)_{p'} \injto \FM^\times$ 
and we denote by 
$\tilde{\imath} : \QM \to \FM^\times$ the composition of the morphisms 
$\QM \longto (\QM/\ZM)_{p'} \longmapright{\imath} \FM^\times$. 
Finally, we set
$\tilde{\imath}_\Tb : \QM \otimes_\ZM Y(\Tb) \longto \Tb$, 
$r \otimes_\ZM \l \longmapsto \l(\tilde{\imath}(r))$. The image of 
$\tilde{\imath}_\Tb$ is the torsion subgroup of $\Tb$. 

To understand the structure of $C_\Gb(s)$, then, by Proposition 
\ref{ordre fini pareil} and by Remark \ref{injectif pi}, 
it is sufficient to work under the following 
hypothesis~:

\bigskip

\begin{quotation}
{\noindent {\bf Hypothesis - } {\it From now on, and until the 
end of this paper, we assume that $\Gb$ is semisimple and that $s$ has finite order.}}
\end{quotation}

\bigskip

\noindent{\sc Remarque - } 
It must be noticed that, in view of classifying quasi-isolated semisimple elements, 
this hypothesis is not restrictive (see Remark \ref{injectif pi} and 
Corollary \ref{ordre fini}).\finl

\bigskip

\soussection{Preliminaries\label{subsection not}} 
Let $V$ be the $\QM$-vector space $\QM \otimes_\ZM Y(\Tb)$ 
and let $V^*$ be its dual, identified with $\QM \otimes_\ZM X(\Tb)$. We denote 
by $<,> : V^* \times V \to \QM$ the canonical perfect pairing between $V$ and $V^*$. 
Then $Y(\Tb_\simp)$ may be identified with $<\Phi^\vee>$ and $X(\Tb_\adj)$ 
may be identified with $<\Phi>$. Since $\Gb$ is semisimple, $\D$ is a basis of $V^*$. 
Let $(\varpi_\a^\vee)_{\a \in \D}$ be its dual basis. Then $Y(\Tb_\adj)$ 
may be identified with $\oplus_{\a \in \D} \ZM \varpi_\a^\vee$. 
As expected, we have $Y(\Tb_\simp) \incl Y(\Tb) \incl Y(\Tb_\adj) \incl V=
\QM \otimes_\ZM Y(\Tb_\simp)$. If $v \in V$, let 
$\t_v : V \to V$, $x \mapsto x+v$ denote the translation by $v$. 

Let us recall the following elementary fact~:

\bigskip

\begin{lem}\label{zzo}
The map $Y(\Tb_\adj) \to \Tb$, $\l \mapsto \tilde{\imath}_\Tb(\l)$ 
induces an  isomorphism $(Y(\Tb_\adj)/Y(\Tb))_{p'} \simeq \Zb(\Gb)$. 
The map $Y(\Tb) \to \Tb_\simp$, $\l \mapsto \tilde{\imath}_{\Tb_\simp}(\l)$ 
induces an isomorphism $(Y(\Tb)/Y(\Tb_\simp))_{p'} \simeq \Ker \pi_\simp$. 
\end{lem}

\bigskip

If $\l \in V$, we set
$$\Phi(\l)=\{\a \in \Phi~|~<\a,\l> \in \ZM\}$$ 
$$W_\Gb(\l)=\{w \in W~|~w(\l)-\l \in Y(\Tb)\}.\leqno{\text{and}}$$
We denote by $o_\simp(\l)$ (respectively $o_\adj(\l)$, respectively $o_\Gb(\l)$) 
the order of the image of $\l$ in $V/Y(\Tb_\simp)$ (respectively $V/Y(\Tb_\adj)$, 
respectively $V/Y(\Tb)$).
Let $W^\circ(\l)$ denote the Weyl group of the closed subsystem $\Phi(\l)$ 
of $\Phi$. Then $W^\circ(\l)$ is a normal subgroup of $W_\Gb(\l)$. If we fix a 
positive root system $\Phi^+(\l)$ in $\Phi(\l)$, then we can define 
$$A_\Gb(\l)=\{w \in W_\Gb(\l)~|~w(\Phi^+(\l))=\Phi^+(\l)\}.$$
Then 
$$W_\Gb(\l)=A_\Gb(\l) \ltimes W^\circ(\l).$$
The next lemma shows that, in order to understand the structure of $C_\Gb(s)$, 
it is necessary and sufficient to understand the structure of $W(\l)$, $W^\circ(\l)$ 
and $A_\Gb(\l)$. 

\bigskip

\begin{lem}\label{p'}
Let $\l \in \ZMP \otimes_\ZM Y(\Tb_\simp) \incl V$ and let $s=\tilde{\imath}_\Tb(\l)$. 
Then 
\begin{itemize}
\itemth{a} $o_\Gb(\l)$ is the order of $s$~;

\itemth{b} $\Phi(\l)=\Phi(s)$ so $W^\circ(\l)=W^\circ(s)$. 

\itemth{c} $W_\Gb(\l)=W(s)$ and, if $\Phi^+(\l)=\Phi^+(s)$, 
then $A_\Gb(\l)=A(s)$.
\end{itemize}
\end{lem}

\bigskip

By analogy, we say that $\l$ is {\it $\Gb$-isolated} (\resp {\it $\Gb$-quasi-isolated}) 
if $W^\circ(\l)$ (\resp $W(\l)$) is not contained in a proper parabolic subgroup 
of $W$.

\medskip

Let $W_{\text{aff}}=W \ltimes Y(\Tb_\simp)$ denote the affine Weyl 
group of $\Phi$. If $\l \in V$, we set
$$W_{\text{aff}}(\l)=\{w \in W_{\text{aff}}~|~w(\l)=\l\}.$$
For the proof of the next proposition, see 
\cite[Lemme 13.14 and Remark 13.15 (i)]{dmbook} and 
\cite[Chapter VI, \SEC{2}, Exercise 1]{bourbaki}.

\bigskip

\begin{prop}\label{caracterisation wzero}
Let $\l \in V$. Then
\begin{itemize}
\itemth{a} $W_{\mathrm{aff}}(\l)$ is generated by affine reflections. 
Its image in $W$ is $W^\circ(\l)$.

\itemth{b} $W^\circ(\l)$ is the kernel of the map $W_\Gb(\l) \to Y(\Tb)/Y(\Tb_\simp)$, 
$w \mapsto w(\l)-\l + Y(\Tb_\simp)$.

\itemth{c} The exponent of $A_\Gb(\l)$ divides $o_\Gb(\l)$.
\end{itemize}
\end{prop}

\bigskip

\noindent{\sc Remark - } By Proposition \ref{ordre fini pareil}, 
by Lemma \ref{p'} and by Proposition \ref{caracterisation wzero} 
we get that the centralizer of a semisimple element in a 
simply connected group is connected (Steinberg's Theorem).\finl

\bigskip

\soussection{Affine Dynkin diagram} 
We recall here some results from \cite[Chapter VI, \SEC 2]{bourbaki} 
concerning the affine Dynkin diagram associated to a root system. 
We denote by $\Phi_1$, $\Phi_2$,\dots, $\Phi_r$ the distinct irreducible 
components of $\Phi$.

\medskip 

Let us fix $i \in \{1,2,\dots,r\}$. Let 
$V_i=\QM  \otimes_\ZM <\Phi_i>$. 
Let $W_i$ denote the Weyl group 
of $\Phi_i$. We set $\D_i=\D \cap \Phi_i$, $\Phi_i^+=\Phi^+ \cap \Phi_i$. 
Then $V_i=\oplus_{\a \in \D_i} \QM \varpi_\a^\vee$. 
We denote by $\alpt_i$ the highest root of $\Phi_i$ (with respect 
to the height defined by $\D_i$). We write
$$\alpt_i=\sum_{\a \in \D_i} n_\a \a,$$
where the $n_\a$ are non-zero natural numbers ($\a \in \D_i$). 
By convention, we set $\varpi^\vee_{-\alpt_i}=0$, $n_{-\alpt_i}=1$. 
Let $\Delt_i=\D \cup \{-\alpt_i\}$, 
$\D_{i,\mini}=\{\a \in \D_i~|~n_\a=1\}$ and 
$\Delt_{i,\mini}=\D_{i,\mini} \cup \{-\alpt_i\}$. 
If $\a \in \Delt_{i,\mini}$, we denote by $\Phi_\a$ the parabolic subsystem of $\Phi_i$ 
with basis $\D_i-\{\a\}$ (for instance, $\Phi_{-\alpt_i}=\Phi_i$) 
and we set $\Phi_\a^+=\Phi_i^+ \cap \Phi_\a$. Let $W_\a$ denote the Weyl 
group of the root system $\Phi_\a$ and $w_\a$ its unique element such that 
$w_\a(\Phi_\a^+)=-\Phi_\a^+$. We set $z_\a=w_\a w_{-\alpt} \in W_i$ 
(note that $z_{-\alpt}=1$) and 
$$\Aut_{W_i}(\Delt_i)=\{z \in W_i~|~z(\Delt_i)=\Delt_i\}.$$
By \cite[chapter VI, \SEC 2, Proposition 6]{bourbaki}, we have 
\equat\label{description autw}
\Aut_{W_i}(\Delt_i)=\{z_\a~|~ \a \in \Delt_{i,\mini}\}.
\endequat
If $\a \in \D_i$, we set $m_\a=0$. We also set $m_{-\alpt_i}=-1$. 
Now, let $\CC_i$ denote the alcove 
\eqna
\CC_i&=&\{\l \in V_i~|~\forall~\a \in \Delt_i,~<\a,\l> \ge m_\a\} \\
&=&\{\l \in V_i~|~(\forall~\a \in \D_i,~<\a,\l> \ge 0){\text{~and~}} 
<\alpt_i,\l> \le 1\}.
\endeqna
Then $\CC_i$ is a fundamental domain for the action of the affine Weyl group 
$W_{i,\text{aff}}= W_i \ltimes <\Phi_i^\vee>$ on $V_i$.  
Moreover, $\CC_i$ is a closed simplex with vertices 
$(\varpi_\a^\vee/n_\a)_{\a \in \Delt_i}$. 

\medskip

With the above notation, we have~:
$$W=W_1 \times W_2 \times \dots \times W_r,$$
$$V=V_1 \oplus V_2 \oplus \dots \oplus V_r,$$
$$Y(\Tb_\simp)=\Oplus_{i=1}^r \Bigl(V_i \cap Y(\Tb_\simp)\Bigr)$$
$$Y(\Tb_\adj)=\Oplus_{i=1}^r \Bigl(V_i \cap Y(\Tb_\adj)\Bigr).\leqno{\text{and}}$$
We set $\Delt=\Delt_1 \cup \Delt_2 \cup \dots \cup \Delt_r$. Now, let 
$$\AC=\{z \in W~|~z(\Delt)=\Delt\}.$$
In other words, $\AC$ is the automorphism group of the affine 
Dynkin diagram of $\Gb$ induced by an element of $W$. We have
$$\AC=\Aut_{W_1}(\Delt_1) \times \Aut_{W_2}(\Delt_2) \times 
\dots \times \Aut_{W_r}(\Delt_r).$$
If $z=(z_{\a_1},z_{\a_2},\dots,z_{\a_r}) \in \AC$, with 
$\a_i \in \Delt_{i,\mini}$, we set 
$$\varpi^\vee(z)=\varpi_{\a_1}^\vee + \varpi_{\a_2^\vee} + \dots + \varpi_{\a_r}^\vee.$$
Finally, let 
\eqna
\CC&=&\{\l \in V~|~\forall~\a \in \Delt,~<\a,\l> \ge m_\a\} \\
&=& \CC_1 \times \CC_2 \times \dots \times \CC_r.
\endeqna
Then $\CC$ is a fundamental domain for the action of 
$W_{\text{aff}}$ in $V$. 
Then, by \cite[Chapter VI, \SEC 2]{bourbaki}, we have, for every 
$z \in \AC$, 
\equat\label{stabilisateur C}
z(\CC)+\varpi^\vee(z)=\CC
\endequat
and the map 
\equat\label{iso centre}
\fonction{\varpi^\vee}{\AC}{Y(\Tb_\adj)/Y(\Tb_\simp)}{z}{\varpi^\vee(z) +
Y(\Tb_\simp)}
\endequat
is an isomorphism of groups. If 
$z=(z_{\a_1},z_{\a_2},\dots,z_{\a_r}) \in \AC$, with 
$\a_i \in \Delt_{i,\mini}$, and if $\a \in \Delt_i$, then 
\equat\label{egalite n}
n_{z(\a)}=n_\a
\endequat
and
\equat\label{transformation sommet}
z(\frac{1}{n_\a} \varpi_\a^\vee)+ \varpi_{\a_i}^\vee = 
\frac{1}{n_\a} \varpi_{z(\a)}^\vee.
\endequat
Since we will be working with the affine Weyl group of $W_{\mathrm{aff}}$, it 
will be convenient to work with ``affine coordinates''. More precisely, 
if $\l \in V$, we will denote by $(\l_\a)_{\a \in \Delt}$ the unique 
family of rational numbers such that 
\begin{quotation}
\noindent (1) $\forall~i \in \{1,2,\dots,r\}$, $\DS{\sum_{\a \in \Delt_i}} \l_\a=1$~;

\noindent (2) $\l=\DS{\sum_{\a \in \Delt} \frac{\l_\a}{n_\a} \varpi_\a^\vee}$.
\end{quotation}
Note that $\l \in \CC$ \ssi $\l_\a \ge 0$ for every $ \a \in \Delt$. 
Then, we have, for every $\a \in \Delt$, 
\equat\label{produit scalaire}
<\a,\l>=\frac{\l_\a}{n_\a}+m_\a.
\endequat

\bigskip

\noindent{\sc Proof of \ref{produit scalaire} - } Recall that 
$m_\a$ has been defined in \SEC\ref{subsection not}. If 
$\a \in \D$, then $m_\a=0$ and, by (3), $<\a,\l>=\l_\a/n_\a$. 
On the other hand, if $\a \in \Delt-\D$, then $m_\a=-1$ and 
there exists a unique 
$i \in \{1,2,\dots,r\}$ such that $\a=-\alpt_i$. Therefore, by (2) and (3), 
$<\a,\l>=-\sum_{\b \in \D_i} \l_\b = \l_\a - 1$.\fin

\bigskip

Moreover, it follows from \ref{transformation sommet} that, 
for every $z \in \AC$,  
\equat\label{action autw}
z(\l)+\varpi^\vee(z)=\sum_{\a \in \Delt} \frac{\l_{z^{-1}(\a)}}{n_\a} \varpi_\a^\vee
\endequat
In other words, $(z(\l)+\varpi^\vee(z))_\a = \l_{z^{-1}(\a)}$. 

\bigskip

\soussection{Orbits under the action of ${\boldsymbol{W \ltimes Y(\Tb)}}$}
Let $\AC_\Gb$ be the subgroup of $\AC$ defined to be 
the inverse image of $Y(\Tb)/Y(\Tb_\simp)$ under the isomorphism 
$\varpi^\vee$. Since $\CC$ is a fundamental domain for the action 
of $W_{\mathrm{aff}}$, it will be interesting to understand whenever 
two elements of $\CC$ are in the same orbit under $W \ltimes Y(\Tb)$. 
The answer is given by the following proposition. 

\bigskip

\begin{prop}\label{cgos autw}
Let $\l$ and $\m$ be two elements of $\CC$ and let $w \in W$. 
If $\m-w(\l) \in Y(\Tb)$, then there exists 
$w^\circ \in W^\circ(\l)$ and 
$z \in \AC_\Gb$ such that $w=zw^\circ$. Moreover, if $d$ is 
a common multiple of $o(\l)$ and $o(\m)$, then $z^d=1$. 
\end{prop}

\bigskip

\proof Assume that $\m-w(\l) \in Y(\Tb)$. Then there exists $z \in \AC_\Gb$ 
and $u \in Y(\Tb_\simp)$ such that $w(\l)-\m = -\varpi^\vee(z) + u$. But, 
$(\t_{-u} w)(\l)=\m-\varpi^\vee(z) \in \CC - \varpi^\vee(z)=z(\CC)$ 
(see \ref{stabilisateur C}). Therefore, $(z^{-1} \t_{-u} w)(\l) \in \CC$. 
Since $\CC$ is a fundamental 
domain for the action of $W_{\text{aff}}$ on $V$ and since 
$z^{-1} \t_{-u} w \in W_{\text{aff}}$, we deduce that $z^{-1} \t_{-u} w(\l)=\l$. 
So, by Proposition \ref{caracterisation wzero} (a), $z^{-1}w \in W^\circ(\l)$, as 
expected. 

For the last assertion, note that the hypothesis implies that 
$d(w(\l)-\m) \in Y(\Tb_\simp)$. Therefore, 
$d \varpi^\vee(z) \in Y(\Tb_\simp)$. Since the map \ref{iso centre} is an 
isomorphism, we get that $z^d=1$.\fin

\bigskip

\begin{coro}\label{conjugue ?}
Let $\l$ and $\m$ be two elements of $\CC$. Then the following assertions are 
equivalent~:
\begin{itemize}
\itemth{1} $\l$ and $\m$ are in the same $W \ltimes Y(\Tb)$-orbit.

\itemth{2} There exists $z \in \AC_\Gb$ such that $z(\l)-\m \in Y(\Tb)$.
\end{itemize}
\end{coro}

\bigskip

\proof Clear.\fin

\bigskip

\soussection{The group ${\boldsymbol{W_\Gb(\l)}}$} 
%
%
Let us now come back to the aim of this section, namely the description 
of the group $W_\Gb(\l)$. Since $\CC$ is a fundamental domain for the action of 
$W_{\text{aff}}$ in $V$, it is sufficient to understand 
the structure of $W_\Gb(\l)$ whenever $\l \in \CC$.

\bigskip

\begin{prop}\label{wgl}
Let $\l \in \CC$. We set 
$I_\l=\{\a \in \Delt~|~\l_\a=0\}=\{\a \in \Delt~|~<\a,\l> = m_\a\}$. 
Then~:
\begin{itemize}
\itemth{a} $I_\l$ is a basis of $\Phi(\l)$. 

\itemth{b} If $\Phi^+(\l)$ is the poitive root system of $\Phi(\l)$ 
associated to the basis $I_\l$, then 
$$A^\Gb(\l)=\{z \in \AC_\Gb~|~\forall~\a \in \Delt,~\l_{z(\a)}=\l_\a\}.$$
\end{itemize}
\end{prop}

\bigskip

\proof (a) For $\a \in \Delt$, let $H_\a=\{v \in V~|~<\a,v> = m_\a\}$. 
Then $(H_\a)_{\a \in \Delt}$ is the family of walls of the alc\^ove 
$\CC$. Moreover, $W_{\mathrm{aff}}$ is generated by the affine reflections 
with respect to the walls of $\CC$ which contains $\l$ (see \cite[??]{bourbaki}). 
Therefore, $W^\circ(\l)$ is generated by the reflections $(s_\a)_{\a \in I_\l}$. 
Since $<\a,\b^\vee> \le 0$ for every $\a$, $\b \in \Delt$, this implies 
that $I_\l$ is a basis of $\Phi(\l)$.

\medskip

(b) Let $A=\{z \in \AC_\Gb~|~\forall~\a \in \Delt,~\l_{z(\a)}=\l_\a\}$. 
Then $A$ stabilizes $I_\l$ by construction and, for every $z \in A$, 
$z(\l)-\l=\varpi^\vee(z) Y(\Tb)$ by \ref{action autw}. So $A \incl A_\Gb(\l)$. 
Let us prove now the reverse inclusion. 

First, let us prove that $A_\Gb(\l) \incl \AC_\Gb$.
Let $z \in A_\Gb(\l)$. By Proposition \ref{cgos autw}, there exists 
$a \in \AC_\Gb$ and $w^\circ \in W^\circ(\l)$ such that $z=aw^\circ$.
So $a \in W(\l)$ and $a(I_\l) \incl \Delt$. In particular, 
$a(I_\l) \incl \Phi(\l) \cap \Delt$. But 
$\Phi(\l) \cap \Delt=I_\l$ by (a). 
So $a(I_\l) = I_\l$. Moreover, $z(I_\l)=I_\l$ by definition of 
$A_\Gb(\l)$. So $w^\circ(I_\l)=I_\l$ and $w^\circ \in W^\circ(\l)$, 
which implies that $w^\circ=1$, that is $z=a$. This shows that 
$z \in \AC_\Gb$. 

Now, by \ref{action autw}, we have 
$$z(\l)-\l + \varpi^\vee(z) = \sum_{i=1}^r\Bigl(
\sum_{\a \in \D_i} \frac{\l_{z^{-1}(\a)}-\l_\a}{n_\a} \varpi_\a^\vee \Bigr) \in 
Y(\Tb) \incl Y(\Tb_\adj).$$
Since $z$ stabilizes $I_\l$, we have, for every $\a \in \D$, 
$$\l_{z^{-1}(\a)}\l_\a=0 \imp \l_{z^{-1}(\a)}=\l_\a=0.$$
Moreover, $0 \le \l_a \le 1$. 
Therefore, $(\l_{z^{-1}(\a)}-\l_\a)/n_\a \in ]-1/n_\a,1/n_\a[$. 
Moreover, $(\varpi_\a^\vee)_{\a \in \D}$ is a $\ZM$-basis of 
$Y(\Tb_\adj)$. So $\l_{z^{-1}(\a)}=\l_\a$ for every $\a \in \D$. 
Then, by condition (1), $\l_{z^{-1}(\a)}=\l_\a$ for every $\a \in \Delt$.\fin

\bigskip

\remark{contre exemple ags autw}
Keep the notation of Proposition \ref{wgl}. Then it may happen 
that $A_\Gb(\l)$ is strictly contained in the stabilizer of $I_\l$ in 
$\AC_\Gb$. Take 
for instance $\Gb=\Pb\Gb\Lb_2(\FM)$ and $\l=\varpi_\a^\vee/3$ where 
$\a$ is the unique simple root of $\Gb$.\finl

\bigskip

\remark{I}
If $\l \in \CC$, note that $I_\l \cap \Delt_i \not= \Delt_i$ for 
every $i \in \{1,2,\dots,r\}$.\finl

\bigskip

If $\l \in \CC$, we will choose for $\Phi^+(\l)$ 
the positive root subsystem of $\Phi(\l)$ associated to the basis $I_\l$.

\bigskip

\sursection{Classification of quasi-isolated elements}

\bigskip

\soussection{A characterization of quasi-isolated elements} 
If $I$ is a subset of $\Delt$ such that $I \cap \Delt_i \not= \Delt_i$ 
for every $i \in \{1,2,\dots,r\}$, we denote by $\Phi_I$ the root subsystem 
of $\Phi$ with basis $I$ and by $W_I$ the Weyl group of $\Phi_I$. 
It must be noticed that $W_I$ is not necessarily a \para of $W$. 
The Proposition \ref{wgl} shows that, whenever $\l \in \CC$, 
$W_\Gb(\l)=A \ltimes W_{I_\l}$ for some subgroup $A$ of $\AC$ 
stabilizing $I_\l$. To determine if such a subgroup is contained 
or not in a proper \para of $W$, we need to determine the dimension 
of its fixed-points space. This is done in general in the next lemma.

\bigskip

\begin{lem}\label{dimension espace fixe}
Let $I$ be a subset of $\Delt$ such that $I \cap \Delt_i \not= \Delt_i$ 
for every $i \in \{1,2,\dots,r\}$ and let 
$A$ be a subgroup of $\AC$ stabilizing $I$. 
Let $r'$ denote the number of orbits of $A$ in $\Delt-I$. Then 
$\dim_\QM V^{A \ltimes W_I} = r'-r$.
\end{lem}

\bigskip

\proof 
By taking direct products, we may assume that $\Phi$ is irreducible or, 
in other words, that $r=1$. Let 
$V_I=\QM \otimes_\ZM <\Phi_I>$ and let $E_I$ be the orthogonal 
of $I$ in  $V$. Then $V=V_I \oplus E_I$ and $AW_I^\circ$ stabilizes 
$V_I$ and $E_I$. Moreover, 
$$\{v \in V_I~|~\forall~w \in W_I,~w(v)=v\}=\{0\}$$
and $W_I$ acts trivially on $E_I$. Consequently,
$$V^{A \ltimes W_I}=E_I^A.$$
Let $\QM[\Delt-I]$ denote the $\QM$-vector space with basis 
$(e_\a)_{\a \in \Delt-I}$. 
This is a permutation $A$-module. Let $f : \QM[\Delt-I] \to E_I$ 
the $\QM$-linear map sending $e_\a$ on the projection of 
$\a$ in $E_I$ (for every $\a \in \Delt-I$). 
Then $f$ is a morphism of $\QM A$-modules, whose kernel 
has dimension $1$ (because $|\Delt|=\dim V+1$). 

Let $M=\{v \in \QM[\Delt-I]~|~\forall z \in A,~z(v)=v\}$. Then 
$$\dim_\QM E_I^A=\dim_\QM M - \dim_\QM(M \cap \Ker f).$$
Since $\dim_\QM M = r'$, we only need to show that $A$ acts trivially 
on $\Ker f$.
But $\sum_{\a \in \Delt} n_\a \a = 0$. So, 
by projection on $E_I$, we get that $\Ker f$ is generated by 
$\sum_{\a \in \Delt - I} n_\a e_\a$. 
By equality \ref{egalite n}, this element is invariant under the action of $A$. 
This completes the proof of Lemma \ref{dimension espace fixe}.\fin

\bigskip

\begin{coro}\label{qisole potentiel}
Let $I$ be a subset of $\Delt$ such that $I \cap \Delt_i \not= \Delt_i$ 
for every $i \in \{1,2,\dots,r\}$ and let $A$ be a subgroup of $\AC$ 
stabilizing $I$. Then $A.W_I$ is not contained in a proper parabolic subgroup 
of $W$ \ssi $A$ acts transitively on $\Delt_i-I$ for every $i \in \{1,2,\dots,r\}$.
\end{coro}

\bigskip

\proof This follows immediately from Proposition \ref{dimension espace fixe}.\fin

\bigskip

\begin{coro}\label{test quasiisole}
Let $\l \in \CC$. Then~:
\begin{itemize}
\itemth{a} $\l$ is $\Gb$-isolated \ssi $|\Delt_i - I_\l|=1$ 
for every $i \in \{1,2,\dots,r\}$. 

\itemth{b} $\l$ is $\Gb$-quasi-isolated \ssi $A_\Gb(\l)$ acts 
transitively on $\Delt_i - I_\l$ for every $i \in \{1,2,\dots,r\}$.
\end{itemize}
\end{coro}

\bigskip

\proof This follows immediately from Proposition \ref{wgl} 
and Corollary \ref{qisole potentiel}.\fin

\bigskip

\soussection{Classification of quasi-isolated elements in ${\boldsymbol{V}}$}
We are now ready to complete the classification of conjugacy classes of 
$\Gb$-quasi-isolated elements in $V$. 
Let $\QC(\Gb)$ denote the set of subsets $\O$ of $\Delt$, 
such that, for every $i \in \{1,2,\dots,r\}$, $\O \cap \Delt_i \not=\vide$ 
and the stabilizer of $\O_i$ in $\AC_\Gb$ acts transitively 
on $\Delt_i$. If $\O$ is such a subset, we set
$$\l_\O=\sum_{i=1}^r \Bigl(\frac{1}{n_i(\O) |\O \cap \Delt_i|} 
\sum_{\a \in \O \cap \Delt_i} \varpi_\a^\vee\Bigr),$$
where $n_i(\O)$ is equal to $n_\a$ for every $\a \in \O \cap \Delt_i$ (see equality 
\ref{egalite n}). Note that $\AC_\Gb$ acts on $\QC(\Gb)$. 
Moreover, by \ref{transformation sommet}, we have, for every $z \in \AC_\Gb$, 
\equat\label{transformation lambda o}
z(\l_\O)+\varpi^\vee(z) = \l_{z(\O)}.
\endequat
Finally, we denote by $o_i^\Gb(\O)$ the number $o_\Gb(\varpi_\a^\vee)$ where 
$\a \in \O \cap \Delt_i$. Note that this number is constant on $\O \cap \Delt_i$. 
All the work done in this section shows that~:

\bigskip

\begin{theo}\label{quasi isole lambda}
With the above notation, we have~:
\begin{itemize}
\itemth{a} The map $\QC(\Gb) \to \CC$, $\O \mapsto \l_\O$ induces a bijection between 
the set of orbits of $\AC_\Gb$ in $\QC(\Gb)$ and the set 
of $W \ltimes Y(\Tb)$-orbits of quasi-isolated elements in $V$. 

\itemth{b} Let $\O \in \QC(\Gb)$. Then~:
\begin{itemize}
\itemth{\a} $W^\circ(\l_\O)=W_{\Delt-\O}$~;

\itemth{\b} $A_\Gb(\l_\O)=\{z \in \AC_\Gb~|~z(\O)=\O$~;

\itemth{\g} $o_\Gb(\l)$ is the lowest common multiple 
of $(n_i(\O) o_i^\Gb(\O) |\O \cap \Delt_i|)_{1 \le i \le r}$~;

\itemth{\d}
$\l_\O$ is $\Gb$-isolated \ssi $|\O_i|=1$ for every $i \in \{1,2,\dots,r\}$. 
\end{itemize}
\end{itemize}
\end{theo}

\bigskip

\soussection{Classification of quasi-isolated semisimple elements} 
Let $\Delt_{p'}$ denote the subset of elements $\a \in \Delt$ such that 
$\varpi_\a^\vee/n_\a \in \ZMP \otimes_\ZM Y(\Tb_\simp)$. 
Let $\QC(\Gb)_{p'}$ denote the set of $\O \in \QC(\Gb)$ such that 
$\O \incl \Delt_{p'}$ and, for every $i \in \{1,2,\dots,r\}$, 
$p$ does not divide $|\O \cap \Delt_i|$. If $\O \in \QC(\Gb)_{p'}$, 
we set $t_\O=\tilde{\imath}_\Tb(\l_\O) \in \Tb$.

\bigskip

\begin{theo}\label{quasi isole s}
With the above notation, we have~:
\begin{itemize}
\itemth{a} The map $\QC(\Gb)_{p'} \to \Tb$, $\O \mapsto t_\O$ 
induces a bijection between the set of orbits of $(\AC_\Gb)_{p'}$ 
in $\QC(\Gb)_{p'}$ and the set 
of conjugacy classes of quasi-isolated semisimple elements in $\Gb$.
 
\itemth{b} If $\O \in \QC(\Gb)_{p'}$ then~:
\begin{itemize}
\itemth{\a} $W^\circ(t_\O)=W_{\Delt-\O}$~;

\itemth{\b} $A_\Gb(t_\O)=\{z \in \AC_\Gb~|~z(\O)=\O\}$~;

\itemth{\g} $o(t_\O)$ is the lowest common multiple 
of $(n_i(\O) o_i^\Gb(\O) |\O \cap \Delt_i|)_{1 \le i \le r}$~;

\itemth{\d}
$t_\O$ is $\Gb$-isolated \ssi $|\O_i|=1$ for every $i \in \{1,2,\dots,r\}$. 
\end{itemize}
\end{itemize}
\end{theo}

\bigskip

\proof By Theorem \ref{quasi isole lambda} and Lemma 
\ref{p'}, it is enough to show that the map $\QC(\Gb)_{p'} \to \CC$, 
$\O \mapsto \l_\O$ induces a bijection between the set 
of orbits of $(\AC_\Gb)_{p'}$ in $\QC(\Gb)_{p'}$ 
to the set of $W \ltimes Y(\Tb)$-orbits of quasi-isolated elements $\l$ in $V$ 
such that $p$ does not divide $o_\simp(\l)$. But this follows from 
Theorem \ref{quasi isole lambda} (b) ($\g$) and the last assertion of 
Proposition \ref{cgos autw}.\fin

\bigskip

\remark{tres bon} We recall that the prime number $p$ is said to be {\it very good} 
for $\Gb$ if it does not divide the numbers $n_\a$ ($\a \in \D$) 
and $|\AC|=|Y(\Tb_\simp)|/|Y(\Tb_\adj)|$. 
We say here that $p$ is {\it almost very good} 
for $\Gb$ is it does not divide the numbers $n_\a$ ($\a \in \Delt$) and 
$|\AC_\Gb|=|Y(\Tb)|/|Y(\Tb_\adj)|$. If $p$ is very good, then it 
is almost very good. 

If $p$ is almost very good, then 
$\Delt_{p'}=\Delt$ and $(\AC_\Gb)_{p'}=\AC_\Gb$ so the set of 
$W \ltimes Y(\Tb)$-orbits of $\Gb$-quasi-isolated elements in $V$ is 
in natural bijection with the set of conjugacy classes of quasi-isolated 
semisimple elements in $\Gb$ (through the map $\tilde{\imath}_\Tb$).\finl

\bigskip

\example{p=2} If all the irreducible components of $\Phi$ are of type 
$B$, $C$ ou $D$ and if $p=2$, then $\Delt_{p'}=\{-\alpt_1,-\alpt_2,\dots,-\alpt_r\}$. 
Therefore, $1$ is the unique 
quasi-isolated element in $\Gb$.\finl

\bigskip

\soussection{Simply connected groups} 
If $\Gb$ is simply connected then $\AC_\Gb=\{1\}$. 
Therefore, we retrieve the well-known classification of isolated semisimple 
elements in $\Gb$~:

\bigskip

\begin{prop}\label{simplement connexe}
Assume that $\Gb$ is semisimple and simply connected. 
Then the map $\Delt_{1,p'} \times \dots \times \Delt_{r,p'} \to \Gb$, 
$(\a_1,\dots,\a_r) \mapsto \prod_{i=1}^r t_{\varpi_{\a_i}^\vee/n_{\a_i}}$ 
induces a bijection between 
$\Delt_{1,p'} \times \dots \times \Delt_{r,p'}$ and the set of conjugacy classes 
of (quasi-)isolated elements in $\Gb$.
\end{prop}

\bigskip

\example{symplectique}
Assume here that $p\not=2$ and that $\Gb=\Sb\pb(V)$ where $V$ is an even-dimensional 
vector space endowed with a non-degenerate alternating form. Let $\dim V=2n$. 
Then $\Gb$ is simply connected, so $\AC_\Gb=\{1\}$. Moreover, 
$\Delt_{p'}=\Delt$. Let us write $\a_0=-\alpt_1$ and let us number the affine 
Dynkin diagram $\Delt$ of $\Gb$ as follows
$$
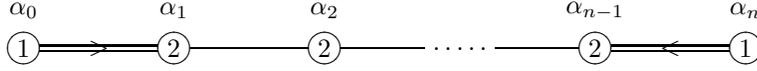
\xy 
(0,0) *++={1} *\frm{o}; (20,0) *++={2} *\frm{o} **@{=}; 
(40,0) *++={2} *\frm{o} **@{-}; (58,0) *++={.\dots.} **@{-}; 
(76,0) *++={2} *\frm{o} **@{-}; (96,0) *++={1} *\frm{o} **@{=};
(0,5) *++={\a_0}; (20,5) *++={\a_1}; 
(40,5) *++={\a_2}; (76,5) *++={\a_{n-1}}; (96,5) *++={\a_n}; 
(10,0) *++={>}; (86,0) *++={<}; 
\endxy$$ 
The natural numbers written inside the node $\a_i$ is the number $n_{\a_i}$. 
For $0 \le i \le n$, let $\O_i=\{\a_i\}$ and let 
$t_i=\tilde{\imath}_\Tb(\varpi_{\a_i}^\vee/n_{\a_i})$. 
Then $\{t_i~|~0 \le i \le n\}$ is a set of representatives of 
conjugacy classes of isolated (i.e. quasi-isolated) elements in $\Gb$. 
Note that $t_i$ is characterized by the following two properties~:
$$t_i^2=1\quad\text{ and } \quad \dim \Ker(t_i+\Id_V)=i.$$
This shows that an element $s$ is isolated in $\Gb$ \ssi $s^2=1$. Finally, 
note that $C_\Gb(t_i)=\Sb\pb_{2i}(\FM) \times \Sb\pb_{2(n-i)}(\FM)$.\fin

\bigskip

\soussection{Special orthogonal groups\label{sous special orthogonal}}
The case of special orthogonal groups in characteristic $2$ has been treated in 
Example \ref{p=2}. In this subsection, we study the case of special 
orthogonal groups in good characteristic. We first adopt a naive point-of-view, 
using the natural representation of special orthogonal groups. At the end 
of this subsection, we will explain the link between this point-of-view 
and Theorem \ref{quasi isole s}. 

\bigskip

\begin{quotation}
\noindent{\bf Hypothesis : } {\it 
Let us assume in this subsection, and only in this subsection, 
that $p\not=2$ and that $\Gb=\Sb\Ob(V,\langle,\rangle)=\Sb\Ob(V)$ where 
$V$ is a finite dimensional vector space over $\FM$ and $\langle,\rangle$ 
is a non degenerate symmetric bilinear form on $V$.}
\end{quotation} 

\bigskip

We denote by $n$ the rank of $\Gb$. Then $n=\Bigl[\DS{\frac{\dim V}{2}}\Bigr]$, 
except whenever $\dim V=2$ (in this case, $n=0$). If $s^2=1$, then 
$\dim \Ker(s+\Id_V) \equiv 0 \mod 2$ so 
$\dim \Ker(s-\Id_V) \equiv \dim V \mod 2$. 

\bigskip

\begin{prop}\label{carre}
With this hypothesis, we have~:
\begin{itemize}
\itemth{a} $s$ is quasi-isolated \ssi $s^2=1$.

\itemth{b} If $s^2=1$, then $s$ is isolated \ssi $\dim \Ker(s - \e \Id_V) \not= 1$ 
for every $\e \in \{1,-1\}$.
\end{itemize}
\end{prop}

\bigskip

\proof Assume first that there exists an eigenvalue $\z$ of $s$ such that 
$\z^2\not= 1$. Let $V_\z$ 
denote the $\z$-eigenspace of $s$ in $V$. Let $E$ be the orthogonal 
subspace to $V_\z \oplus V_{\z^{-1}}$. We have 
$$V=V_\z \oplus V_{\z^{-1}} \oplus E,$$
and this is an orthogonal decomposition. Therefore, the centralizer of $s$ 
in $\Gb$ is contained in 
$\Gb \cap (\Gb\Lb(V_\z) \times \Gb\Lb(V_{\z^{-1}}) \times \Gb\Lb(E) )$, 
which is a \levi of a proper \para of $\Gb$. So $s$ is not 
quasi-isolated. 

\medskip

Assume now that $s^2=1$. Then $V=V_1 \oplus V_{-1}$ 
and this decomposition is orthogonal. So 
\equat
C_\Gb(s)=(\Ob(V_1) \times \Ob(V_{-1})) \cap \Gb\quad\text{et}\quad 
C_\Gb^\circ(s)=\Sb\Ob(V_1) \times \Sb\Ob(V_{-1}).
\equat
So $s$ is quasi-isolated and it is isolated \ssi 
$\dim V_1 \not=1$ and $\dim V_{-1} \not= 1$.\fin

\bigskip

\begin{coro}\label{special classe}
Keep the hypothesis of this subsection. 
If $0 \le i \le n$, let $t_i$ denote a semisimple element of $\Gb$ 
such that $t_i^2=1$ and $\dim \Ker(t_i+\Id_V) = 2i$. 
Then $\{t_i~|~0 \le i \le n\}$ 
is a set of representatives of conjugacy classes of quasi-isolated elements in 
$\Gb$. Moreover, $t_i$ is isolated \ssi $i \not\in \{1, (\dim V)/2 -1\}$.
\end{coro}

\bigskip

Let us now compare the description given by Corollary \ref{special classe} 
and the one given by Theorem \ref{quasi isole s}. Since $p \not= 2$, 
we have $\Delt_{p'}=\Delt$ and $\AC_{p'}=\AC$ 
(indeed, $p$ is very good for $\Gb$). For getting a uniform description, 
we assume that $\dim V \not\in\{1,2,3,4,6\}$ (whenever $\dim V \in \{1,2,3,4,6\}$, 
then the reader can also check that Corollary \ref{special classe} and Theorem 
\ref{quasi isole s} are still compatible~!). 

\bigskip

\soussoussection{Type $B$} 
We assume here that $\dim V=2n+1$ and that $n \ge 2$. 
We set $\a_0=-\alpt_1$. Then $\AC_\Gb=\AC$ is of order $2$. 
We denote by $\s$ its unique non-trivial element. 
We number the affine Dynkin diagram of $\Gb$ as follows~: 
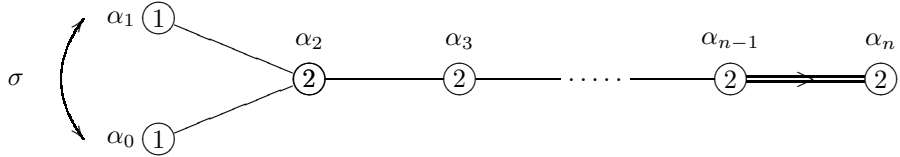
$$\xy 
(0,8) *++={1} *\frm{o} ; (20,0) *++={2} *\frm{o} **@{-};
(0,-8) *++={1} *\frm{o} ; (20,0) *++={2} *\frm{o} **@{-};
(40,0) *++={2} *\frm{o} **@{-}; (58,0) *++={.\dots.} **@{-}; 
(76,0) *++={2} *\frm{o} **@{-}; (96,0) *++={2} *\frm{o} **@{=};
(-5,8) *++={\a_1}; (-5,-8) *++={\a_0}; 
(-10,8); (-10,-8) **\crv{(-16,0)} ?<*\dir{<} ?>*\dir{>}; 
(-19,0) *++={\s}; (20,5) *++={\a_2}; 
(40,5) *++={\a_3}; (76,5) *++={\a_{n-1}}; (96,5) *++={\a_n}; (86,0) *++={>};
\endxy$$ 
The natural number written inside the node $\a_i$ is equal to $n_{\a_i}$. 
We have 
$$\s=\lexp{s_{\a_1} s_{\a_2} \dots s_{\a_{n-1}}}{s_{\a_n}}.$$ 
The action of $\s$ on $\Delt$ is given by the above 
diagram~: $\s(\a_0)=\a_1$, $\s(\a_1)=\a_0$, 
$\s(\a_i)=\a_i$ for every $i \in \{2,3,\dots, n\}$.

\medskip

Let $\O_0=\{\a_0\}$, $\O_1=\{\a_0,\a_1\}$ and, for $2 \le i \le n$, 
let $\O_i=\{\a_i\}$. Then $\O_i \in \QC(\Gb)$. Moreover, one can check that 
$\{\O_0,\O_1,\dots,\O_n\}$ is a set of representatives of 
$\AC_\Gb$-orbits in $\QC(\Gb)$. 
Then
$$t_{\O_i}^2=1 \quad\text{ and } \quad \dim \Ker(t_{\O_i}+\Id_V)=i.$$
This shows that $t_{\O_i}=t_i$~: we retrieve Corollary \ref{special classe}.

\bigskip

\soussoussection{Type $D$} 
We assume here that $\dim V=2n$ and that $n \ge 4$. 
We set $\a_0=-\alpt_1$. Then $\AC_\Gb$ is of order $2$. 
We denote by $\s$ its unique non-trivial element. Note that 
$\AC_\Gb \not= \AC$. 
We number the affine Dynkin diagram of $\Gb$ as follows~: 
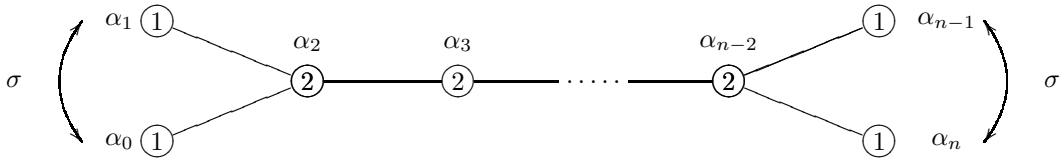
$$\xy 
(0,8) *++={1} *\frm{o} ; (20,0) *++={2} *\frm{o} **@{-};
(0,-8) *++={1} *\frm{o} ; (20,0) *++={2} *\frm{o} **@{-};
(40,0) *++={2} *\frm{o} **@{-}; (58,0) *++={.\dots.} **@{-}; 
(76,0) *++={2} *\frm{o} **@{-}; (96,8) *++={1} *\frm{o} **@{-};
(76,0) *++={2} *\frm{o} **@{-}; (96,-8) *++={1} *\frm{o} **@{-};
(-5,8) *++={\a_1}; (-5,-8) *++={\a_0}; 
(-10,8); (-10,-8) **\crv{(-16,0)} ?<*\dir{<} ?>*\dir{>}; 
(-19,0) *++={\s};
(20,5) *++={\a_2}; (40,5) *++={\a_3}; (76,5) *++={\a_{n-2}}; 
(105,8) *++={\a_{n-1}}; (105,-8) *++={\a_n}; 
(110,8) ; (110,-8) **\crv{(116,0)} ?<*\dir{<} ?>*\dir{>} ;
(119,0) *++={\s};
\endxy$$
The natural number written inside the node $\a_i$ is equal to $n_{\a_i}$. We have 
$$\s=\lexp{s_{\a_1} s_{\a_2} \dots s_{\a_{n-2}}}{(s_{\a_{n-1}}s_{\a_n})}.$$ 
The action of $\s$ on $\Delt$ is given in the above diagram~: 
$\s(\a_0)=\a_1$, that $\s(\a_1)=\a_0$ 
and that $\s(\a_i)=\a_i$ for every $i \in \{2,3,\dots,n-2\}$, $\s(\a_{n-1})=\a_n$ 
and $\s(\a_n)=\a_{n-1}$. 
 
\medskip

Let $\O_0=\{\a_0\}$, $\O_1=\{\a_0,\a_1\}$, 
$\O_i=\{\a_i\}$ (for $2 \le i \le n-2$), $\O_{n-1}=\{\a_{n-1},\a_n\}$ 
and $\O_n=\{\a_n\}$. Then $\O_i \in \QC(\Gb)$. Moreover, one can check that 
$\{\O_0,\O_1,\dots,\O_n\}$ is a set of representatives of 
$\AC_\Gb$-orbits in $\QC(\Gb)$.
Then
$$t_{\O_i}^2=1 \quad\text{ and } \quad \dim \Ker(t_{\O_i}+\Id_V)=i.$$
This shows that $t_{\O_i}=t_i$~: we retrieve Corollary \ref{special classe}.

\bigskip

\sursection{Adjoint simple groups}

\bigskip

The aim of this section is to provide complete tables for isolated 
and quasi-isolated semisimple conjugacy classes in adjoint simple groups. 
For classical groups, we also give a description in terms of their natural 
representation.

\bigskip

\begin{quotation}
\noindent{\bf Hypothesis :} {\it In this section, and only in this section, 
we assume that $\Gb$ is adjoint and simple.} 
\end{quotation}

\bigskip

The root system $\Phi$ is then irreducible. We denote by $\a_0$ the root $-\alpt_1$. 
Note also that $\AC_\Gb=\AC$. Moreover, for 
every $\a \in \Delt$, we have $o_\Gb(\varpi_\a^\vee)=1$. Therefore, 
the Theorem \ref{quasi isole s} can be stated as follows~:

\bigskip

\begin{theo}\label{quasi isole adjoint}
Assume that $\Gb$ is adjoint and simple. Then $\QC(\Gb)_{p'}$ is the set 
of subsets $\O$ in $\Delt_{p'}$ which are acted on transitively by 
their stabilizer in $\AC$. If $\O \in \QC(\Gb)_{p'}$, let $n_\O$ 
denote the number $n_\a$ (for some $\a \in \O$). Then~:
$$t_\O=
\tilde{\imath}_\Tb(\frac{1}{n_\O.|\O|} \sum_{\a \in \O} \varpi_\a^\vee).$$ 
We have~:
\begin{itemize}
\itemth{a} The map $\QC(\Gb)_{p'} \to \Tb$, $\O \mapsto t_\O$ 
induces a bijection between the set of orbits of $\AC_{p'}$ 
in $\QC(\Gb)_{p'}$ and the set 
of conjugacy classes of quasi-isolated semisimple elements in $\Gb$.
 
\itemth{b} Let $\O \in \QC(\Gb)_{p'}$. Then~:
\begin{itemize}
\itemth{\a} $W^\circ(t_\O)=W_{\Delt-\O}$~;

\itemth{\b} $A_\Gb(t_\O)=\{z \in \AC_\Gb~|~z(\O)=\O\}$~;

\itemth{\g} $o(t_\O)=n_\O|\O|$~;

\itemth{\d}
$s_\O$ is $\Gb$-isolated \ssi $|\O|=1$. 
\end{itemize}
\end{itemize}
\end{theo}

\bigskip

This implies that the set of conjugacy classes of isolated semisimple 
elements in $\Gb$ is in bijection with the set of orbits of $\AC$ 
in $\Delt_{p'}$. 

\medskip

\soussection{Classification by use of the affine Dynkin diagram} 
We first set some notation. We denote by $\a_0$ the root $\-\alpt_1$ 
(recall that $r=1$). Let $n$ denote the rank of $\Gb$ (i.e. $n=|\D|$). 
We write $\D=\{\a_1,\a_2,\dots,\a_n\}$. If $0 \le i \le n$, we set 
$n_{\a_i}=n_i$, $z_{\a_i}=z_i$ and $\varpi_{\a_i}^\vee=\varpi_i^\vee$. 
Note that $z_i(\a_0)=\a_i$. 
The Table I gives the list of all the affine Dynkin diagrams 
together with the structure of $\AC$ (see \cite[Planches I-IX]{bourbaki}).

\medskip

We give in Table II the classification of conjugacy classes of quasi-isolated 
elements in adjoint classical groups. In Table III, we deal with the 
adjoint groups of exceptional type $E_6$ and $E_7$. We have not included adjoint 
groups of type $E_8$, $F_4$ and $G_2$ since they are also simply connected. 
Therefore, Proposition \ref{simplement connexe}, Theorem 
\ref{quasi isole adjoint} and Table I gives easily all informations 
concerning the (quasi-)isolated elements for these groups.

\medskip

\soussection{Explicit descriptions for adjoint classical groups} 
The case of special orthogonal groups was done in subsection 
\ref{sous special orthogonal}. Therefore, we only have to investigate 
adjoint classical groups of type $A$, $C$ and $D$.

\medskip

\soussoussection{Type $A$}
Assume here that $\Gbt=\Gb\Lb_{n+1}(\FM)$, that $\Gb=\Pb\Gb\Lb_{n+1}(\FM)$ and that 
$\pi : \Gbt \to \Gb$ is the canonical morphism (here, $n$ is a non-zero natural 
number). Let $I_{n+1}$ denote the identity matrix of 
$\Gb\Lb_{n+1}(\FM)$. If $d$ is a non-zero natural number invertible in $\FM$, 
we denote by $\z_d$ a primitive $d$-th root of unity in $\FM^\times$ and we 
set $J_d=\diag(1,\z_d,\z_d^2,\dots,\z_d^{d-1}) \in \Gb\Lb_d(\FM)$. 

Now, let $\Div_{p'}(n+1)$ denote the set of divisors of $n+1$ which are invertible 
in $\FM$. If $d \in \Div_{p'}(n+1)$, let $\sti_{n+1,d}$ denote the matrix 
$I_{\frac{n+1}{d}} \otimes J_d \in \Gbt$. We set $s_{n+1,d}=\pi(\sti_{n+1,d})$. 

Note that $\AC_{p'}$ is cyclic of order $n_{p'}=|\Delt_{p'}|$ 
and that it acts transitively on $\Delt_{p'}$. If $d \in \Div_{p'}(n+1)$, 
we denote by $\O_{n+1,d}$ the orbit of $\a_0$ under the unique subgroup 
of order $d$ of $\AC$~: $\O_{n+1,d}=\{\a_{j(n+1)/d}~|~0 \le j \le d-1\}$.  

\medskip

\begin{prop}\label{pgl}
If $\Gb=\Pb\Gb\Lb_{n+1}(\FM)$, then the map $\Div_{p'}(n+1) \to \Gb$, 
$d \mapsto s_{n+1,d}$ 
is a bijection between $\Div_{p'}(n+1)$ and the set of conjugacy classes of 
quasi-isolated semisimple elements in $\Gb$. Through the parametrization 
of Theorem \ref{quasi isole adjoint}, this corresponds to the map 
$\Div_{p'}(n+1) \mapsto \QC(\Gb)_{p'}$, $d \mapsto \O_{n+1,d}$. 

If $d \in \Div_{p'}(n+1)$, then $s_{n+1,d}$ has order $d$, 
$W^\circ(s)\simeq (\SG_{n+1/d})^d$ 
and $A(s) \simeq (\ZM/d\ZM)$ acts on $W^\circ(s)$ by permutation of the 
components. Moreover, $s_{n+1,d}$ is isolated if and only if $d = 1$.
\end{prop}

\medskip

\soussoussection{Type $C$} 
We assume here that $p \not= 2$. 
Let $V$ be a $2n$-dimensional vector space over $\FM$, with $n \ge 2$. 
Let $\b : V \times V \to \FM$ be a non-degenerate skew-symmetric 
bilinear form. We assume here that $\Gbt=\Sb\pb(V,\b)$, 
that $\Gb=\Gbt/\{\Id_V,-\Id_V\}$ and that $\pi : \Gbt \to \Gb$ 
is the canonical morphism. 

\medskip

\begin{prop}\label{psp}
Let $\sti \in \Gbt$ be semisimple and let $s=\pi(\sti)$. 
\begin{itemize}
\itemth{a} If $s$ is quasi-isolated, then $\sti^4=1$.

\itemth{b} If $\sti^2=1$, then $s$ is isolated. 

\itemth{c} If $\sti^4=1$ and $\sti^2\not=1$, then 
$s$ is quasi-isolated if and only if 
$\dim \Ker(\sti - \Id_V)=\dim \Ker(\sti+\Id_V)$.
\end{itemize}
\end{prop}

\medskip 

\proof (a) follows immediately from Corollary \ref{quasi isole} and 
Example \ref{symplectique}. (b) and (c) 
follow from direct computation.\fin

\medskip

If $0 \le i \le n/2$, let $\tti_i$ be an element of $\Gbt$ 
such that $\tti^2_i=1$ and $\dim \Ker(\tti_i + \Id_V) = 2i$. 
If $0 \le i < n/2$, let $\sti_i$ be an element of $\Gbt$ 
of order $4$ such that $\dim \Ker(\sti_i-\Id_V)=\dim \Ker(\sti_i+\Id_V)=i$. 
We set $t_i=\pi(\tti_i)$ and $s_i=\pi(\sti_i)$. 

\medskip

\begin{coro}\label{classification psp}
The set $\{t_i~|~0 \le i \le n/2\} \cup \{s_i~|~0 \le i < n/2\}$ is 
a set of representatives of quasi-isolated elements of $\Gb$. 
The subset of $\Delt$ associated to $t_i$ 
(respectively $s_i$) through the parametrization of Theorem 
\ref{quasi isole adjoint} is $\{\a_i\}$ (respectively $\{\a_i,\a_{n-i}\}$).
\end{coro}

\medskip

\soussoussection{Type $D$} 
We assume here that $p \not= 2$. 
Let $V$ be a $2n$-dimensional vector space over $\FM$, with $n \ge 3$. 
Let $\b : V \times V \to \FM$ be a non-degenerate symmetric 
bilinear form. We assume here that $\Gbt=\Sb\Ob(V,\b)$, 
that $\Gb=\Gbt/\{\Id_V,-\Id_V\}$ and that $\pi : \Gbt \to \Gb$ 
is the canonical morphism. 

\medskip

\begin{prop}\label{pso}
Let $\sti \in \Gbt$ be semisimple and let $s=\pi(\sti)$. 
\begin{itemize}
\itemth{a} If $s$ is quasi-isolated, then $\sti^4=1$.

\itemth{b} If $\sti^2=1$, then $s$ is quasi-isolated. Moreover, 
$s$ is isolated if and only if $\dim \Ker(\sti-\Id_V) \not\in \{1,n-1\}$.

\itemth{c} If $\sti^4=1$ and $\sti^2\not=1$, then 
$s$ is quasi-isolated if and only if 
$\dim \Ker(\sti - \Id_V)=\dim \Ker(\sti+\Id_V)$ and 
$\dim \Ker(\sti-\Id_V) \not= 0$ if $n$ is odd.
\end{itemize}
\end{prop}

\medskip 

\proof (a) follows immediately from Corollary \ref{quasi isole} and 
Proposition \ref{pso}. (b) and (c) follow from direct computations.\fin

\medskip

If $0 \le i \le n/2$, let $\tti_i$ be an element of $\Gbt$ 
such that $\tti^2_i=1$ and $\dim \Ker(\tti_i + \Id_V) = 2i$. 
If $1 \le i < n/2$, let $\sti_i$ be an element of $\Gbt$ 
of order $4$ such that $\dim \Ker(\sti_i-\Id_V)=\dim \Ker(\sti_i+\Id_V)=i$. 
We set $t_i=\pi(\tti_i)$ and $s_i=\pi(\sti_i)$. 
Finally, there are two conjugacy classes of elements $\sti$ of order $4$ 
such that $\dim \Ker(\sti-\Id_V)=\dim \Ker(\sti+\Id_V)=0$ (these 
two conjugacy classes are in correspondence through the non trivial 
automorphism of the Dynkin diagram of $\Gb$)~: we denote by $\sti_0$ and 
$\sti_0'$ some representatives of these two classes. We set $s_0=\pi(\sti_0)$ 
and $s_0'=\pi(\sti_0')$. We set $E_n=\vide$ if $n$ is odd and 
$E_n=\{s_0,s_0'\}$ if $n$ is even. 

\medskip

\begin{coro}\label{classification pso}
The set $\{t_0\} \cup \{t_i~|~2 \le i \le n/2\}$ is 
a set of representatives of isolated elements of $\Gb$. 
The subset of $\Delt$ associated to $t_i$ 
through the parametrization of Theorem 
\ref{quasi isole adjoint} is $\{\a_i\}$.

The set $\{t_1\} \cup \{s_i~|~1 \le i < n/2\} \cup E_n$ is 
a set of representatives of quasi-isolated but non-isolated elements 
of $\Gb$. Through the parametrization of Theorem 
\ref{quasi isole adjoint}, $t_1$ is associated to $\{\a_0,\a_1\}$, 
$s_1$ is associated to $\{\a_0,\a_1,\a_{n-1},\a_n\}$ and 
$s_i$ is associated to $\{\a_i,\a_{n-i}\}$. If $n$ is even, 
$s_0$ and $s_0'$ correspond to $\{\a_0,\a_{n-1}\}$ and 
$\{\a_0,\a_n\}$ (or conversely).
\end{coro}

$$\begin{array}{|c|c|c|c|}
\hline
\vertical \text{Type of } \Gb & \Delt & \AC & |\AC| \\
\hline
\vertical A_n & \xy 
(0,0)   *++={1} *\frm{o}       ; (10,0)  *++={1} *\frm{o} **@{-};
(20,0)  *++={1} *\frm{o} **@{-}; (30,0)  *++={\dots}      **@{-}; 
(40,0)  *++={1} *\frm{o} **@{-}; (50,0)  *++={1} *\frm{o} **@{-};
(0,0)   *++={1} *\frm{o}       ;  (25,-12) *++={1} *\frm{o} **@{-};
(50,0)  *++={1} *\frm{o}       ;  (25,-12) *++={1} *\frm{o} **@{-};
(0,4)   *++={\a_1}             ; (25,-16) *++={\a_0}; 
(10,4)  *++={\a_2}             ; (20,4)  *++={\a_3}; 
(40,4)  *++={\a_{n-1}}         ; (50,4)  *++={\a_n}; 
\endxy
 & <z_1> & n+1 \\
\hline
\vertical B_n & \xy 
(0,4)   *++={1} *\frm{o}       ; (10,0) *++={2} *\frm{o} **@{-};
(0,-4)  *++={1} *\frm{o}       ; (10,0) *++={2} *\frm{o} **@{-};
(20,0)  *++={2} *\frm{o} **@{-}; (30,0) *++={\dots}      **@{-}; 
(40,0)  *++={2} *\frm{o} **@{-}; (50,0) *++={2} *\frm{o} **@{=};
(-5,4)  *++={\a_1}             ; (-5,-4) *++={\a_0}; 
(10,4)  *++={\a_2}             ; (20,4)  *++={\a_3}; 
(40,4)  *++={\a_{n-1}}         ; (50,4)  *++={\a_n}; 
(45,0)  *++={>};
\endxy & <z_1> & 2 \\
\hline
\vertical C_n & 
\xy 
(0,0)  *++={1} *\frm{o}       ; (10,0) *++={2} *\frm{o} **@{=}; 
(20,0) *++={2} *\frm{o} **@{-}; (30,0) *++={\dots}      **@{-}; 
(40,0) *++={2} *\frm{o} **@{-}; (50,0) *++={2} *\frm{o} **@{-};
(60,0) *++={1} *\frm{o} **@{=};
(0,4)  *++={\a_0}             ; (10,4) *++={\a_1}; 
(20,4) *++={\a_2}             ; (40,4) *++={\a_{n-2}}; 
(50,4) *++={\a_{n-1}}         ; (60,4) *++={\a_n}; 
(5,0)  *++={>}                ; (55,0) *++={<}; 
\endxy & <z_n> & 2 \\
\hline
\vertical \xy (0,3) *++={D_n} ; (0,-3) *++={n\text{ even}};\endxy & \xy 
(0,4)  *++={1} *\frm{o}       ; (10,0) *++={2} *\frm{o} **@{-};
(0,-4) *++={1} *\frm{o}       ; (10,0) *++={2} *\frm{o} **@{-};
(20,0) *++={2} *\frm{o} **@{-}; (30,0) *++={\dots} **@{-}; 
(40,0) *++={2} *\frm{o} **@{-};
(50,0) *++={2} *\frm{o} **@{-}; (60,4) *++={1} *\frm{o} **@{-};
(50,0) *++={2} *\frm{o} **@{-}; (60,-4) *++={1} *\frm{o} **@{-};
(-5,4)  *++={\a_1}    ; (-5,-4) *++={\a_0}; 
(10,4)  *++={\a_2}    ; (20,4)  *++={\a_3}; 
(40,4)  *++={\a_{n-3}}; (50,4)  *++={\a_{n-2}};
(67,4)  *++={\a_{n-1}}; (65,-4) *++={\a_n}; 
\endxy & <z_1> \times <z_n> & 4 \\
\hline
\vertical \xy (0,3) *++={D_n} ; (0,-3) *++={n\text{ odd}};\endxy & \xy 
(0,4)  *++={1} *\frm{o}       ; (10,0) *++={2} *\frm{o} **@{-};
(0,-4) *++={1} *\frm{o}       ; (10,0) *++={2} *\frm{o} **@{-};
(20,0) *++={2} *\frm{o} **@{-}; (30,0) *++={\dots} **@{-}; 
(40,0) *++={2} *\frm{o} **@{-};
(50,0) *++={2} *\frm{o} **@{-}; (60,4) *++={1} *\frm{o} **@{-};
(50,0) *++={2} *\frm{o} **@{-}; (60,-4) *++={1} *\frm{o} **@{-};
(-5,4)  *++={\a_1}    ; (-5,-4) *++={\a_0}; 
(10,4)  *++={\a_2}    ; (20,4)  *++={\a_3}; 
(40,4)  *++={\a_{n-3}}; (50,4)  *++={\a_{n-2}};
(67,4)  *++={\a_{n-1}}; (65,-4) *++={\a_n}; 
\endxy & <z_n> & 4 \\
\hline
\vertical E_6 & \xy (0,0) *++={1} *\frm{o} ; (10,0) *++={2} *\frm{o} **@{-};
(20,0) *++={3} *\frm{o} **@{-}; (30,0) *++={2} *\frm{o} **@{-}; 
(40,0) *++={1} *\frm{o} **@{-};
(20,0) *++={3} *\frm{o}; (20,-8) *++={2} *\frm{o} **@{-}; 
(20,-16) *++={1} *\frm{o} **@{-}; 
(0,4) *+={\a_1}; (10,4) *+={\a_3}; (20,4) *+={\a_4}; 
(30,4) *+={\a_5}; (40,4) *+={\a_6}; 
(15,-8) *+={\a_2}; (15,-16) *+={\a_0};
\endxy  & <z_1> & 3 \\
\hline 
\vertical E_7 & \xy 
(0,0)  *++={1} *\frm{o}       ; (10,0) *++={2} *\frm{o} **@{-};
(20,0) *++={3} *\frm{o} **@{-}; (30,0) *++={4} *\frm{o} **@{-}; 
(40,0) *++={3} *\frm{o} **@{-}; (50,0) *++={2} *\frm{o} **@{-}; 
(60,0) *++={1} *\frm{o} **@{-}; 
(30,0) *++={4} *\frm{o}; (30,-10) *++={2} *\frm{o} **@{-}; 
(0,4) *+={\a_0}   ; (10,4) *+={\a_1}  ; (20,4) *+={\a_3}; 
(30,4) *+={\a_4}  ; (40,4) *+={\a_5}  ; (50,4) *+={\a_6}; 
(60,4) *+={\a_7}; 
(25,-10) *+={\a_2}; 
\endxy & <z_7> & 2 \\
\hline 
\vertical E_8 & \xy 
(0,0)  *++={2} *\frm{o}       ; (10,0) *++={4} *\frm{o} **@{-};
(20,0) *++={6} *\frm{o} **@{-}; (30,0) *++={5} *\frm{o} **@{-}; 
(40,0) *++={4} *\frm{o} **@{-}; (50,0) *++={3} *\frm{o} **@{-}; 
(60,0) *++={2} *\frm{o} **@{-}; (70,0) *++={1} *\frm{o} **@{-}; 
(20,0) *++={6} *\frm{o}; (20,-10) *++={3} *\frm{o} **@{-}; 
(0,4)  *+={\a_1}  ; (10,4) *+={\a_3}  ; (20,4) *+={\a_4}; 
(30,4) *+={\a_5}  ; (40,4) *+={\a_6}  ; (50,4) *+={\a_7}; 
(60,4) *+={\a_8}  ; (70,4) *+={\a_0}  ; (15,-10) *+={\a_2}; 
\endxy & 1 & 1 \\
\hline 
\vertical F_4 & \xy 
(0,0)  *++={1} *\frm{o}       ; (10,0) *++={2} *\frm{o} **@{-};
(20,0) *++={3} *\frm{o} **@{-}; (30,0) *++={4} *\frm{o} **@{=}; 
(40,0) *++={2} *\frm{o} **@{-}; 
(0,4)  *+={\a_0}; (10,4) *+={\a_1}; (20,4) *+={\a_2}; 
(30,4) *+={\a_3}; (40,4) *+={\a_4}; (25,0) *++={>};
\endxy & 1 & 1 \\
\hline 
\vertical G_2 & \xy 
(0,0)  *++={1} *\frm{o}       ; (12,0) *++={2} *\frm{o} **@{-};
(24,0) *++={3} *\frm{o} **@{=};
(0,4)  *+={\a_0}; (12,4) *+={\a_1}; (24,4) *+={\a_2}; 
(18,0) *++={>};
\endxy & 1 & 1 \\
\hline 
\end{array}$$

\bigskip

\begin{centerline}{\sc Table I. 
Affine Dynkin diagrams}\end{centerline}

\newpage

$$\begin{array}{|c|l|c|c|c|c|c|}
\hline
\vertical \quad\Gb\quad & \hskip1cm\O & p~? & 
o(s_\O) & C_\Gb^\circ(s_\O) & 
|A(s_\O)| & \text{isolated~?} \\
\hline
\vertical A_n & \{\a_{j(n+1)/d}~|~0 \le j \le d-1\}
& p \not|~ d & d & (A_{(n+1)/d}-1)^d & d & \text{iff } d=1 \\
& \text{for } d~|~n+1 &&&&& \\
\hline
\vertical B_n & \{\a_0\}  &  & 1 & B_n & 1 & \text{yes} \\
\vertical & \{\a_0,\a_1\} & p\not= 2 & 2 & B_{n-1} & 2 & \text{no} \\
\vertical & \{\a_d\},~2 \le d \le n & p\not= 2 & 2 & D_d \times B_{n-d} & 2 & 
\text{yes} \\
\hline
\vertical C_n & \{\a_0\}            & & 1 & C_n & 1 & \text{yes} \\
\vertical & \{\a_d\},~1 \le d < n/2 & p\not= 2 & 2 & C_d \times C_{n-d} & 
1 & \text{yes} \\
\vertical & \{\a_{n/2}\} \text{ (if $n$ is even)} & 
p\not= 2 & 2 & C_{n/2} \times C_{n/2} & 2 & \text{yes} \\
\vertical &\{\a_0,\a_n\} & p\not= 2 & 2 & A_{n-1} & 2 & \text{no} \\
\vertical &\{\a_d,\a_{n-d}\},~1 \le d < n/2 & 
p\not= 2 & 4 & (B_d)^2 \times A_{n-2d-1} & 2 & \text{no} \\
\hline
\vertical D_n & \{\a_0\}  & & 1 & D_n & 1 & \text{yes} \\
\vertical & \{\a_d\},~2 \le d < n/2 & p\not= 2 & 2 & D_d \times D_{n-d} & 
2 & \text{yes} \\
\vertical & \{\a_{n/2}\} \text{ (if $n$ is even)} & 
p\not= 2 & 4 & D_{n/2} \times D_{n/2} & 4 & \text{yes} \\
\vertical &\{\a_d,\a_{n-d}\},~2 \le d < n/2 & 
p\not= 2 & 4 & (D_d)^2 \times A_{n-2d-1} & 4 & \text{no} \\
\vertical & \{\a_0,\a_1,\a_{n-1},\a_n\} & p\not= 2 & 4 & 
A_{n-3} & 4 & \text{no} \\
\vertical & \{\a_0,\a_1\} & p \not= 2 & 2 & D_{n-1} & 2 & \text{no} \\
\vertical & \{\a_0,\a_{n-1}\} \text{ (if $n$ is even)} 
& p \not = 2 & 2 & A_{n-1} & 2 & \text{no} \\
\vertical & \{\a_0,\a_n\} \text{ (if $n$ is even)} 
& p \not = 2 & 2 & A_{n-1} & 2 & \text{no} \\
\hline
\end{array}$$

\bigskip

\begin{centerline}{\sc Table II. 
Quasi-isolated elements in adjoint classical groups}\end{centerline}

\newpage

$$\begin{array}{|c|c|c|c|c|c|c|}
\hline
\vertical \quad\Gb\quad & \O & p~? & 
o(s_\O) & C_\Gb^\circ(s_\O) & 
|A(s_\O)| & \text{isolated~?} \\
\hline
\vertical E_6 & \{\a_0\} & & 1 & E_6 & 1 & \text{yes} \\
\vertical & \{\a_2\} & p\not= 2 & 2 & 
            A_5 \times A_1 & 1 & \text{yes} \\
\vertical & \{\a_4\} & p\not= 3 & 3 & 
            A_2 \times A_2 \times A_2 & 3 & \text{yes} \\
\vertical & \{\a_0,\a_1,\a_6\} & p \not= 3 & 3 & 
            D_4 & 3 & \text{no} \\
\vertical & \{\a_2,\a_3,\a_5\} & p \not\in \{2,3\} & 6 & 
            A_1 \times A_1 \times A_1 \times A_1 & 3 & \text{no} \\
\hline
\vertical E_7 & \{\a_0\} & & 1 & E_7 & 1 & \text{yes} \\
\vertical & \{\a_1\} & p\not= 2 & 2 & 
            A_1 \times D_6 & 1 & \text{yes} \\
\vertical & \{\a_2\} & p\not= 2 & 2 & 
            A_7 & 2 & \text{yes} \\
\vertical & \{\a_3\} & p\not= 3 & 3 & 
            A_2 \times A_5 & 1 & \text{yes} \\
\vertical & \{\a_4\} & p\not= 2 & 4 & 
            A_3 \times A_3 \times A_1 & 2 & \text{yes} \\
\vertical & \{\a_0,\a_7\} & p \not=2 & 2 & 
            E_6 & 2 & \text{no} \\
\vertical & \{\a_1,\a_6\} & p \not=2 & 4 & 
            D_4 \times A_1 \times A_1 & 2 & \text{no} \\
\vertical & \{\a_3,\a_5\} & p \not\in \{2,3\} & 6 & 
            A_2 \times A_2 \times A_2 & 2 & \text{no} \\
\hline	       	    	       	    
\end{array}$$

\bigskip

\begin{centerline}{\sc Table III. 
Quasi-isolated elements in adjoint groups of type $E_6$ and $E_7$}
\end{centerline}

\end{document}